\documentclass{amsart}
\usepackage{amssymb}

\usepackage{amscd}

\newtheorem{theorem}{Theorem}[section]

\newtheorem{lemma}[theorem]{Lemma}
\newtheorem{proposition}[theorem]{Proposition}

\theoremstyle{definition}
\newtheorem{definition}{Definition}[section]

\theoremstyle{remark}

\begin{document}
\title{Cohomology and Obstructions II: \\
Curves on $K$-trivial threefolds}
\author{Herb Clemens}
\address{Mathematics Department, University of Utah}
\email{clemens@math.utah.edu}
\date{October, 2002}
\maketitle

\begin{abstract}
On a threefold with trivial canonical bundle, Kuranishi theory gives an
algebro-geometry construction of the (local analytic) Hilbert scheme of
curves at a smooth holomorphic curve as a gradient scheme, that is, the
zero-scheme of the exterior derivative of a holomorphic function on a
(finite-dimensional) polydisk. (The corresponding fact in an infinite
dimensional setting was long ago discovered by physicists.) An analogous
algebro-geometric construction for the holomorphic Chern-Simons functional
is presented giving the local analytic moduli scheme of a vector bundle. An
analogous gradient scheme construction for Brill-Noether loci on ample
divisors is also given. Finally, using a structure theorem of
Donagi-Markman, we present a new formulation of the Abel-Jacobi mapping into
the intermediate Jacobian of a threefold with trivial canonical bundle.
\end{abstract}

\section{Introduction\label{intro}\protect\footnote{%
Author partially supported by NSF grant DMS-9970412 }}

\subsection{The problem}

This paper computes the local analytic deformation theory of three closely
related objects on a $K$-trivial threefold $X_{0}$. The three are:

1) Deformations of pairs $\left( X_{0},Y_{0}\right) $ where $Y_{0}$ is a
smooth curve in $X_{0}$.

2) Deformations of pairs $\left( X_{0},E_{0}\right) $ where $E_{0}$ is a
holomorphic vector bundle on $X_{0}$.

3) Deformations of triples $\left( X_{0},S_{0},L_{0}\right) $ where $S$ is a
smooth very ample divisor on $X_{0}$ and $L$ is a line bundle on $S_{0}$
whose Chern class is an algebraic $1$-cycle on $S_{0}$ which goes to zero in 
$H_{2}\left( X_{0};\Bbb{Z}\right) $. In this case we make the additional
assumption that 
\begin{equation*}
h^{1}\left( \mathcal{O}_{X_{0}}\right) =h^{2}\left( \mathcal{O}%
_{X_{0}}\right) =0.
\end{equation*}

The special properties of all these local analytic deformation schemes in
the case of $K$-trivial threefolds derive from the (Serre) duality between
the first-order deformations for \textit{fixed} $K$-trivial threefold $X_{0}$
, given by 
\begin{equation*}
Ext^{1}\left( A,A\right) ,
\end{equation*}
and the obstruction space to extension to higher orders, given by 
\begin{equation*}
Ext^{2}\left( A,A\right) .
\end{equation*}
This duality is given by the natural pairing 
\begin{equation*}
Ext^{1}\left( A,A\right) \times Ext^{2}\left( A,A\right) \rightarrow
Ext^{3}\left( A,A\right)
\end{equation*}
coupled with a trace map 
\begin{equation*}
Ext^{3}\left( A,A\right) \rightarrow H^{3}\left( \mathcal{O}_{X_{0}}\right) =%
\Bbb{C}.
\end{equation*}
This point of view is that of, for example, \cite{KS}. In case 1), for the
inclusion map 
\begin{equation*}
i:Y_{0}\rightarrow X_{0}
\end{equation*}
apply 
\begin{equation*}
R\frak{\hom }_{\mathcal{O}_{X_{0}}}\left( \ ,i_{*}\mathcal{O}_{Y_{0}}\right)
\end{equation*}
to the exact sequence 
\begin{equation*}
0\rightarrow \mathcal{I}_{Y_{0}}\rightarrow \mathcal{O}_{X_{0}}\rightarrow
i_{*}\mathcal{O}_{Y_{0}}\rightarrow 0
\end{equation*}
to obtain 
\begin{equation*}
\frak{ext}^{1}\left( i_{*}\mathcal{O}_{Y_{0}},i_{*}\mathcal{O}
_{Y_{0}}\right) =N_{Y_{0}\backslash X_{0}}.
\end{equation*}
Apply 
\begin{equation*}
R\frak{\hom }_{\mathcal{O}_{X_{0}}}\left( i_{*}\mathcal{O}_{Y_{0}},\ \right)
\end{equation*}
to the exact sequence 
\begin{equation*}
0\rightarrow \mathcal{I}_{Y_{0}}\rightarrow \mathcal{O}_{X_{0}}\rightarrow
i_{*}\mathcal{O}_{Y_{0}}\rightarrow 0
\end{equation*}
to obtain a surjection 
\begin{equation*}
\omega _{Y_{0}}\rightarrow \frak{ext}^{2}\left( i_{*}\mathcal{O}
_{Y_{0}},i_{*}\mathcal{O}_{Y_{0}}\right) .
\end{equation*}
Since a Koszul resolution shows that $\frak{ext}^{2}\left( i_{*}\mathcal{O}%
_{Y_{0}},i_{*}\mathcal{O}_{Y_{0}}\right) $ is locally of rank $1$ we
conclude 
\begin{equation*}
\frak{ext}^{2}\left( i_{*}\mathcal{O}_{Y_{0}},i_{*}\mathcal{O}
_{Y_{0}}\right) =\omega _{Y_{0}}.
\end{equation*}
So 
\begin{equation*}
A=i_{*}\mathcal{O}_{Y_{0}}.
\end{equation*}
In case 2), 
\begin{equation*}
A=E_{0}
\end{equation*}
and the trace map 
\begin{equation*}
H^{3}\left( End\left( E_{0}\right) \right) \rightarrow H^{3}\left( \mathcal{O%
}_{X_{0}}\right)
\end{equation*}
is the obvious one.

Finally, in case 3), let 
\begin{equation*}
j:S_{0}\rightarrow X_{0}
\end{equation*}
denote the inclusion and put 
\begin{equation*}
A=j_{*}L_{0}.
\end{equation*}
Now 
\begin{equation*}
R\frak{\hom }_{\mathcal{O}_{X_{0}}}\left( j_{*}L_{0},j_{*}L_{0}\right) =R%
\frak{\hom }_{\mathcal{O}_{X_{0}}}\left( j_{*}\mathcal{O}_{S},j_{*}\mathcal{O%
}_{S}\right)
\end{equation*}
and, as above, 
\begin{equation*}
\frak{ext}^{1}\left( j_{*}\mathcal{O}_{S},j_{*}\mathcal{O}_{S}\right)
=N_{S\backslash X_{0}}.
\end{equation*}
Also 
\begin{equation*}
\frak{\hom }\left( j_{*}L,j_{*}L\right) =j_{*}\mathcal{O}_{S}.
\end{equation*}
The behavior of $L_{0}$ is reflected at the $E_{2}$-term of the
local-to-global spectral sequence, namely 
\begin{eqnarray*}
E_{2}^{0,1} &=&H^{0}\left( \frak{ext}^{1}\left( j_{*}L_{0},j_{*}L_{0}\right)
\right) =H^{0}\left( N_{S\backslash X_{0}}\right) \\
E_{2}^{2,0} &=&H^{2}\left( \frak{\hom }\left( j_{*}L_{0},j_{*}L_{0}\right)
\right) =H^{2}\left( \mathcal{O}_{S}\right) ,
\end{eqnarray*}
and, letting $\nabla $ denote the Gauss-Manin connection, 
\begin{equation*}
\begin{array}{r}
d_{2}:H^{0}\left( N_{S\backslash X_{0}}\right) \rightarrow H^{2}\left( 
\mathcal{O}_{S}\right) \\ 
\zeta \mapsto \nabla _{\zeta }\left( c_{1}\left( L_{0}\right) \right)
\end{array}
.
\end{equation*}
For more details of this computation (applied to a special case), see the
Appendix of \cite{KS}.

\subsection{The setting for curves on $K$-trivial threefolds}

Let 
\begin{equation*}
X_{0}
\end{equation*}
be a smooth $K$-trivial K\"{a}hler threefold. The deformations of $X_{0}$
are unobstructed (see, for example, \cite{C}). We let 
\begin{equation}
s:X\rightarrow X^{\prime }  \label{a1}
\end{equation}
be a versal deformation of $X_{0}$ over an analytic polydisk $X^{\prime }$,
that is, the natural map 
\begin{equation*}
\left. T_{X^{\prime }}\right| _{0}\hookrightarrow H^{1}\left(
T_{X_{0}}\right)
\end{equation*}
is an isomorphism. Let 
\begin{equation*}
F^{j}H^{3}=H^{3,0}\left( X/X^{\prime }\right) \oplus \ldots \oplus
H^{j,3-j}\left( X/X^{\prime }\right) \subseteq H^{3}\left( X/X^{\prime };%
\Bbb{C}\right)
\end{equation*}
denote the Hodge filtration and let 
\begin{equation*}
n^{\prime }:=\dim X^{\prime }+1=h^{0}\left( \Omega _{X_{0}}^{3}\right)
+h^{1}\left( \Omega _{X_{0}}^{2}\right) .
\end{equation*}

Let $Y_{0}\subseteq X_{0}$ be a smooth irreducible curve which is\textit{\
homologically ample}. By this last we mean that there is a family $%
L/X^{\prime }$ of smooth curves in $X/X^{\prime }$ such that $L_{0}$ is
disjoint from $Y_{0}$ and 
\begin{equation*}
\left\{ Y_{0}\right\} \equiv \left\{ rL_{0}\right\} \in H_{2}\left( X_{0};%
\Bbb{Z}\right)
\end{equation*}
for some $r\in \Bbb{Q}$. We let 
\begin{equation*}
Y^{\prime }
\end{equation*}
denote an analytic neighborhood of $\left\{ Y_{0}\right\} $ in the relative
Hilbert scheme of (proper) curves in $X/X^{\prime }$. Let 
\begin{equation*}
p:Y\rightarrow Y^{\prime }
\end{equation*}
denote the universal curve and 
\begin{equation*}
\pi :Y^{\prime }\rightarrow X^{\prime }
\end{equation*}
the induced map. Then, since the tangent space to the deformation space of
the pair $\left( X_{0},Y_{0}\right) $ is 
\begin{equation*}
\Bbb{H}^{1}\left( T_{X_{0}}\rightarrow N_{Y_{0}\backslash X_{0}}\right) ,
\end{equation*}
which sits in the exact sequence 
\begin{equation}
0\rightarrow H^{0}\left( N_{Y_{0}\backslash X_{0}}\right) \rightarrow \Bbb{H}
^{1}\left( T_{X_{0}}\rightarrow N_{Y_{0}\backslash X_{0}}\right) \rightarrow
H^{1}\left( T_{X_{0}}\right) ,  \label{edge}
\end{equation}
$Y^{\prime }$ can be realized as a closed analytic subscheme of an analytic
polydisk 
\begin{equation*}
U^{\prime }
\end{equation*}
of dimension equal to 
\begin{equation*}
n^{\prime }+h^{0}\left( N_{Y_{0}\backslash X_{0}}\right)
\end{equation*}
for which there is a smooth (surjective) morphism 
\begin{equation*}
U^{\prime }\rightarrow X^{\prime }
\end{equation*}
extending $\pi $ above. Abusing notation we denote this extension again as $%
\pi $. Thus we have an exact sequence 
\begin{equation*}
0\rightarrow T_{\pi }\rightarrow T_{U^{\prime }}\rightarrow \pi
^{*}T_{X^{\prime }}\rightarrow 0.
\end{equation*}
By $\left( \ref{edge}\right) $ we have a natural inclusion of sequences 
\begin{equation}
\begin{array}{ccccccccc}
0 & \rightarrow & H^{0}\left( N_{Y_{0}\backslash X_{0}}\right) & \rightarrow
& \Bbb{H}^{1}\left( T_{X_{0}}\rightarrow N_{Y_{0}\backslash X_{0}}\right) & 
\rightarrow & H^{1}\left( T_{X_{0}}\right) &  &  \\ 
&  & \downarrow ^{\lambda } &  & \downarrow &  & \downarrow ^{\mu } &  &  \\ 
0 & \rightarrow & \left. T_{\pi }\right| _{\left( \left\{ Y_{0}\right\}
,\left\{ X_{0}\right\} \right) } & \rightarrow & \left. \pi ^{*}T_{U^{\prime
}}\right| _{\left( \left\{ Y_{0}\right\} ,\left\{ X_{0}\right\} \right) } & 
\rightarrow & \left. \pi ^{*}T_{X^{\prime }}\right| _{\left( \left\{
Y_{0}\right\} ,\left\{ X_{0}\right\} \right) } & \rightarrow & 
\end{array}
\label{diag}
\end{equation}
such that $\lambda $ and $\mu $ are both isomorphisms.

\subsection{Choice of a $3$-form}

Finally we let 
\begin{equation*}
\tilde{X}^{\prime }
\end{equation*}
denote the analytic manifold obtained by removing the zero-section from the
total space on the analytic line bundle 
\begin{equation*}
s_{*}\Omega _{X/X^{\prime }}^{3}
\end{equation*}
on $X^{\prime }$ and use `tilde' to denote base extension by $\tilde{X}
^{\prime }/X^{\prime }$ , that is, 
\begin{eqnarray*}
\tilde{X} &=&X\times _{X^{\prime }}\tilde{X}^{\prime } \\
\tilde{Y}^{\prime } &=&Y^{\prime }\times _{X^{\prime }}\tilde{X}^{\prime } \\
\tilde{Y} &=&Y\times _{X^{\prime }}\tilde{X}^{\prime } \\
\tilde{\pi } &:&\tilde{U}^{\prime }\rightarrow \tilde{X}^{\prime }.
\end{eqnarray*}

\subsection{The ``potential function'' $\Phi $}

Our first main goal in this paper is to construct (shrinking $X^{\prime }$
and $U^{\prime }$ as necessary) a holomorphic function 
\begin{equation*}
\Phi
\end{equation*}
on $\tilde{U}^{\prime }$ such that, with respect to the exact sequence, 
\begin{equation*}
0\rightarrow \tilde{\pi }^{*}\Omega _{\tilde{X}^{\prime }}^{1}\rightarrow
\Omega _{\tilde{U}^{\prime }}^{1}\rightarrow \Omega _{\tilde{U}^{\prime }/%
\tilde{X}^{\prime }}^{1}\rightarrow 0,
\end{equation*}
we have:

\begin{description}
\item[Property 1]  The relative Hilbert scheme $\tilde{Y}^{\prime }$,
considered as an analytic subscheme of $\tilde{U}^{\prime }$ is the
zero-scheme of the section 
\begin{equation*}
d_{\tilde{U}^{\prime }/\tilde{X}^{\prime }}\Phi
\end{equation*}
of 
\begin{equation*}
\Omega _{\tilde{U}^{\prime }/\tilde{X}^{\prime }}^{1}.
\end{equation*}

\item[Property 2]  Under a natural isomorphism 
\begin{eqnarray*}
F^{2}H^{3}\left( \tilde{X}/\tilde{X}^{\prime }\right) &\cong &T_{\tilde{X}%
^{\prime }} \\
\Omega _{\tilde{X}^{\prime }}^{1} &\cong &\left( F^{2}H^{3}\left( \tilde{X}/%
\tilde{X}^{\prime }\right) \right) ^{\vee }
\end{eqnarray*}
given by Donagi-Markman (see \S \ref{taut} below or \S 1 of \cite{V}), the
section 
\begin{equation*}
\left. d\Phi \right| _{\tilde{Y}^{\prime }}
\end{equation*}
of 
\begin{equation*}
\tilde{\pi }^{*}\Omega _{\tilde{X}^{\prime }}^{1}
\end{equation*}
is the normal function 
\begin{equation}
\int\nolimits_{r\left( L\times _{X^{\prime }}Y^{\prime }\right) /Y^{\prime
}}^{Y/Y^{\prime }}:\tilde{Y}^{\prime }\rightarrow \left( F^{2}H^{3}\left( 
\tilde{X}/\tilde{X}^{\prime }\right) \right) ^{\vee }.  \label{nlfn}
\end{equation}
\end{description}

An immediate corollary of Property 2 is that the image of $\left. d\Phi
\right| _{\tilde{Y}^{\prime }}$ under the natural map 
\begin{equation*}
\left( F^{2}H^{3}\left( \tilde{U}^{\prime }\times _{X^{\prime }}X/\tilde{U}%
^{\prime }\right) \right) ^{\vee }\rightarrow \left( F^{2}H^{3}\left( \tilde{
X}/\tilde{X}^{\prime }\right) \right) ^{\vee }
\end{equation*}
is Lagrangian with respect to the symplectic structure of Donagi-Markman
(again see below or \S 1 of \cite{V}).

\subsection{Holomorphic Chern-Simons and Brill-Noether theory}

Reacting to a preliminary version of the above result for deformations of $%
\left( Y_{0},X_{0}\right) $, both Richard Thomas and Claire Voisin saw wider
settings in which analogous results were true. As Thomas pointed out, an
analogous theorem, proved below, must hold for deformations of $\left(
E_{0},X_{0}\right) $ where the ``holomorphic Chern-Simons functional'' plays
the role of the potential function. Additionally, at the author's
invitation, Richard Thomas authored an additional chapter for this paper
explaining the direct link between the deformation problem for $\left(
Y_{0},X_{0}\right) $ and that for $\left( E_{0},X_{0}\right) $. This is
accomplished through an analogue of Abel's theorem, as anticipated in \cite
{T}. Finally Claire Voisin authored an final chapter establishing the
analogous results for deformations of $\left( X_{0},S,L\right) $, the
Brill-Noether loci for hypersurface sections of a Calabi-Yau threefold.

\subsection{Acknowledgements}

The author wishes to thank Sheldon Katz for suggesting the problem of
realizing Hilbert schemes as gradient schemes, and for providing the
contextual overview presented in the first paragraph of the introduction.
Thanks also to Claire Voisin for an absolutely critical remark that led to
the solution of the problem, as well as for the chapter giving the extension
to Brill-Noether theory. Finally, the author also wishes to thank Richard
Thomas for pointing out that there should be an analogous construction in
holomorphic Chern-Simons theory (see \cite{DT}, \cite{W}, and \S 7 of \cite
{T}) which we treat towards the end of this paper, and for his chapter
explaining the generalization of Abel's theorem relating the holomorphic
Chern-Simons functional to the potential function introduced above.

The author also wishes to thank the Institute for Advanced Study, Princeton,
NJ, for its hospitality during the period in which this paper was written.

\section{The trivialization $F$}

\begin{proposition}
\label{goodprop}Let $\left( 0,0\right) \in U^{\prime }$ represent the
basepoint $\left( \left\{ Y_{0}\right\} ,\left\{ X_{0}\right\} \right) $.
Shrinking $U^{\prime }$ as necessary, there is a $C^{\infty }$-isomorphism 
\begin{equation}
F:\left( U^{\prime }\times _{X^{\prime }}X\right) /U^{\prime }\overset{%
\left( \sigma ,ident._{U^{\prime }}\right) }{\longrightarrow }X_{0}\times
U^{\prime }  \label{triv}
\end{equation}
such that

i) 
\begin{equation*}
\left. F\right| _{\pi ^{-1}\left\{ \left( 0,0\right) \right\}
}=identity_{X_{0}},
\end{equation*}

ii) 
\begin{equation*}
\sigma ^{-1}\left( x_{0}\right)
\end{equation*}
is a analytic submanifold for all $x_{0}\in X_{0}$,

iii) 
\begin{equation*}
\sigma ^{-1}\left( Y_{0}\right) \supseteq Y
\end{equation*}

iv) 
\begin{equation*}
\sigma ^{-1}\left( L_{0}\right) =L\times _{X^{\prime }}U^{\prime }.
\end{equation*}
\end{proposition}

\begin{proof}
The proof is just as in Theorem 13.1 of \cite{C}, using a patching argument
on local data to construct a convergent $C^{\infty }$-isomorphism 
\begin{equation*}
\left( U^{\prime }\times _{X^{\prime }}X\right) /U^{\prime }\rightarrow
X_{0}\times U^{\prime },
\end{equation*}

satisfying i)-iv).
\end{proof}

As in Part One of \cite{C}, the deformation of $\overline{\partial }$
determined by the trivialization $F$ is given by an operator 
\begin{equation*}
\overline{\partial }-L_{\xi }
\end{equation*}
on $X_{0}\times U^{\prime }$ where 
\begin{equation*}
\xi \in A_{X_{0}}^{0,1}\left( T_{X_{0}}\right) \otimes \Bbb{C}\left[ \left[
U^{\prime }\right] \right] .
\end{equation*}
In fact, by Appendix A of \cite{C}, the trivialization has the additional
property that $\xi $ is convergent, that is, $\xi $ is given by a
well-defined holomorphic mapping 
\begin{equation*}
\xi :U^{\prime }\rightarrow A_{X_{0}}^{0,1}\left( T_{X_{0}}^{1,0}\right)
\end{equation*}
comprised of pointwise holomorphic mappings 
\begin{equation*}
\xi \left( x\right) :U^{\prime }\rightarrow \left. \left( T_{X_{0}}^{\vee
}\right) ^{0,1}\right| _{x}\otimes \left. \left( T_{X_{0}}^{1,0}\right)
\right| _{x}
\end{equation*}
for each $x\in X_{0}$. By construction, 
\begin{eqnarray*}
\left. \xi \right| _{Y_{0}\times Y^{\prime }} &\in &A_{Y_{0}}^{0,1}\left(
T_{Y_{0}}\right) \otimes \mathcal{O}_{Y^{\prime }} \\
\left. \xi \right| _{L_{0}\times U^{\prime }} &\in &A_{L_{0}}^{0,1}\left(
T_{L_{0}}\right) \otimes \mathcal{O}_{Y^{\prime }}.
\end{eqnarray*}
An important point to notice is that, if $\mathcal{I}_{Y^{\prime }}$ denotes
the ideal of $Y^{\prime }$ in $U^{\prime }$ and $\frak{m}$ is the (maximal)
ideal of $\left( 0,0\right) \in U^{\prime }$, then 
\begin{equation*}
\left. \xi \right| _{Y_{0}\times U^{\prime }}\in A_{Y_{0}}^{0,1}\left(
\left. T_{X_{0}}\right| _{Y_{0}}\right) \otimes \mathcal{O}_{Y^{\prime }}
\end{equation*}
gives the element in 
\begin{equation*}
H^{1}\left( N_{Y_{0}\backslash X_{0}}\right) \otimes \frac{\mathcal{I}
_{Y^{\prime }}}{\frak{m}\cdot \mathcal{I}_{Y^{\prime }}}
\end{equation*}
which is the obstruction to extending the family $Y/Y^{\prime }$. (See \cite
{C}, \S 11.) Furthermore $F$ determines a $C^{\infty }$-mapping 
\begin{equation*}
f:Y_{0}\times U^{\prime }\rightarrow X\times _{X^{\prime }}U^{\prime }
\end{equation*}
by the rule 
\begin{equation*}
f\left( u^{\prime }\right) =\left. F^{-1}\right| _{Y_{0}\times \left\{
u^{\prime }\right\} }.
\end{equation*}
We denote the family of $C^{\infty }$-deformations of $Y_{0}$ as 
\begin{equation*}
U/U^{\prime }
\end{equation*}
with fibers $U_{u^{\prime }}$ which are no longer algebraic curves but
simply smooth $2$-real-dimensional manifolds when $u^{\prime }\notin
Y^{\prime }$.

\section{Hodge theory in terms of the trivialization $F$}

\subsection{The holomorphic $3$-form on $X/X^{\prime }$}

Let $\eta $ be an everywhere non-zero global section of 
\begin{equation*}
F^{3}H^{3}\left( X/X^{\prime }\right) .
\end{equation*}
Abusing notation we shall contiue to denote by $\eta $ the pullback of this
form to 
\begin{equation*}
F^{3}H^{3}\left( \left( U^{\prime }\times _{X^{\prime }}X\right) /U^{\prime
}\right) .
\end{equation*}
In Appendix A of \cite{C} the trivialization 
\begin{equation*}
F:\left( U^{\prime }\times _{X^{\prime }}X\right) /U^{\prime }\rightarrow
X_{0}\times U^{\prime }
\end{equation*}
in Proposition \ref{goodprop} is constructed the property that 
\begin{equation*}
\left( F^{-1}\right) ^{*}\left( \eta \right)
\end{equation*}
is a \textit{holomorphic} family of $3$-forms on $X_{0}$, that is, is given
by a holomorphic mapping from $U^{\prime }$ into the (infinite-dimensional)
complex vector space 
\begin{equation*}
A_{X_{0}}^{3}
\end{equation*}
of global $C^{\infty }$-three-forms on $X_{0}$. Thus, if $\Gamma $ is a
fixed $3$-chain in $X_{0}$, then by differentiating under the integral sign, 
\begin{equation*}
\int\nolimits_{\Gamma }\left( F^{-1}\right) ^{*}\left( \eta \right)
\end{equation*}
is a holomorphic function on $U^{\prime }$.

\subsection{Abel-Jacobi map}

As in \cite{C}, we define for the trivialization given by Proposition \ref
{goodprop} the ``Hodge spaces'' 
\begin{equation*}
K^{p,q}=H^{q}\left( A_{X_{0}}^{p,*},\overline{\partial }_{X_{0}}-L_{\xi
}^{1,0}\right) .
\end{equation*}
Then in \cite{C}, Lemma 8.2, it is shown that the $C^{\infty }$-isomorphism $%
F$ induces the formal correspondence 
\begin{equation*}
H^{q}\left( \Omega _{\left( X\times _{X^{\prime }}U^{\prime }\right)
/U^{\prime }}^{p}\right) \leftrightarrow e^{\left\langle \left. \xi \right|
\ \right\rangle }K^{p,q}.
\end{equation*}
Here it is important to note that this correspondence is to be understood
modulo $\overline{u^{\prime }}$ where $u^{\prime }$ is a system of
holomorphic coordinates on $U^{\prime }$ centered at our given basepoint $%
\left( \left\{ Y_{0}\right\} ,\left\{ X_{0}\right\} \right) $--indeed, for $%
q>0$, the Hodge spaces $H^{q}\left( \Omega _{\left( X\times _{X^{\prime
}}U^{\prime }\right) /U^{\prime }}^{p}\right) $ do not vary holomorphically
in $u^{\prime }$. So writing 
\begin{equation*}
\eta _{0}=\left. \eta \right| _{X_{0}}
\end{equation*}
there is an element 
\begin{equation*}
\beta \in A_{X_{0}}^{3,0}\otimes \frak{m}\cdot \Bbb{C}\left[ \left[
U^{\prime }\right] \right]
\end{equation*}
such that 
\begin{equation}
\eta _{0}+\beta \in K^{3,0}  \label{beta}
\end{equation}
and, modulo $\overline{u^{\prime }}$, 
\begin{equation*}
\eta =F^{*}\left( e^{\left\langle \left. \xi \right| \ \right\rangle }\left(
\eta _{0}+\beta \right) \right) .
\end{equation*}

Also we can rewrite the pull-back of the (relative) intermediate Jacobian 
\begin{equation*}
\mathcal{J}\left( X/X^{\prime }\right) =\frac{\left( F^{2}H^{3}\left(
X/X^{\prime }\right) \right) ^{\vee }}{H_{3}\left( X/X^{\prime };\ \Bbb{Z}
\right) }
\end{equation*}
to a complex torus bundle over $U^{\prime }$ as 
\begin{equation*}
\mathcal{J}\left( \left( X\times _{X^{\prime }}U^{\prime }\right) /U^{\prime
}\right) =\frac{Hom_{\mathcal{O}_{X^{\prime }}}\left( e^{\left\langle \left.
\xi \right| \ \right\rangle }\left( K^{3,0}\oplus K^{2,1}\right) ,\mathcal{O}
_{X^{\prime }}\right) }{H_{3}\left( X_{0};\Bbb{Z}\right) }.
\end{equation*}
In this last formulation the normal function associated to the cycle 
\begin{equation*}
Y/Y^{\prime }-L\times _{X^{\prime }}Y^{\prime }
\end{equation*}
is given by the restriction to $Y^{\prime }\subseteq \left( U^{\prime
}\times X^{\prime }\right) $ of the (formal) function 
\begin{equation}
\int\nolimits_{rL_{0}}^{Y_{0}}:e^{\left\langle \left. \xi \right| \
\right\rangle }\left( K^{3,0}\oplus K^{2,1}\right) \rightarrow \Bbb{C}.
\label{goodmap}
\end{equation}
Let $\Gamma $ be a rational $3$-chain such that 
\begin{equation*}
\partial \Gamma =Y_{0}-rL_{0}.
\end{equation*}
So for our choice of $\eta $ above the holomorphic function $%
\int\nolimits_{\Gamma }\left( F^{-1}\right) ^{*}\eta $ can be rewritten as 
\begin{equation}
\int\nolimits_{\Gamma }\left( e^{\left\langle \left. \xi \right| \
\right\rangle }\left( \eta _{0}+\beta \right) \right) .  \label{bestyet}
\end{equation}

\section{The `potential function' $\Phi $ for curves on $K$-trivial
threefolds}

\subsection{The tautological section\label{taut}}

Following \cite{DM} (see \S 1 of \cite{V}), $F^{3}H^{3}\left( \tilde{X}/%
\tilde{X}^{\prime }\right) $ has a tautologial section 
\begin{equation*}
\tau :\left( x^{\prime },\omega _{x^{\prime }}\right) \mapsto \omega
_{x^{\prime }}
\end{equation*}
and, following Donagi-Markman \cite{DM}, the Gauss-Manin connection $\nabla $
induces a holomorphic bundle isomorphism 
\begin{equation}
\begin{array}{c}
\nabla :T_{\tilde{X}^{\prime }}\rightarrow F^{2}H^{3}\left( \tilde{X}/\tilde{%
X}^{\prime }\right) \\ 
\zeta \mapsto \nabla _{\zeta }\tau
\end{array}
\label{Gauss}
\end{equation}
which restricts to an isomorphism 
\begin{equation*}
T_{\tilde{X}^{\prime }/X^{\prime }}\rightarrow F^{3}H^{3}\left( \tilde{X}/%
\tilde{X}^{\prime }\right) .
\end{equation*}
Since $\nabla $ is a real operator, we extend $\left( \ref{Gauss}\right) $
to a real isomorphism 
\begin{equation}
T_{\tilde{X}^{\prime }}\oplus \overline{T_{\tilde{X}^{\prime }}}\rightarrow
F^{2}H^{3}\left( \tilde{X}/\tilde{X}^{\prime }\right) \oplus \overline{
F^{2}H^{3}\left( \tilde{X}/\tilde{X}^{\prime }\right) }=H^{3}\left( \tilde{X}
/\tilde{X}^{\prime };\Bbb{C}\right)  \label{realGauss}
\end{equation}
which restricts to 
\begin{equation*}
T_{\tilde{X}^{\prime }}\left( \Bbb{R}\right) \rightarrow H^{3}\left( \tilde{X%
}/\tilde{X}^{\prime };\Bbb{R}\right) .
\end{equation*}

\subsection{Distinguished coordinates for $\tilde{X}^{\prime }$}

The isomorphism $\left( \ref{Gauss}\right) $ implies that, for our non-zero
holomorphic section $\eta _{0}$ of $F^{3}H^{3}\left( X_{0}\right) $, the
maps 
\begin{equation}
\tilde{x}_{i}=\left( \int\nolimits_{\gamma _{i}}\tau -\int\nolimits_{\gamma
_{i}}\eta _{0}\right)  \label{goodcoords}
\end{equation}
give distinguished local coordinates for $\tilde{X}^{\prime }$ centered at $%
\left( 0,\eta _{0}\right) $. If 
\begin{equation*}
\left\{ \omega _{i}=pr_{F^{2}H^{3}\left( \tilde{X}/\tilde{X}^{\prime
}\right) }\left( \gamma ^{i}\right) \right\} _{i=1,\ldots ,n^{\prime }}
\end{equation*}
denotes the framing of $F^{2}H^{3}\left( \tilde{X}/\tilde{X}^{\prime
}\right) $ dual to $\left\{ \gamma _{i}\right\} _{i=1,\ldots ,n^{\prime }}$,
then we have 
\begin{equation*}
\omega _{i}=\frac{\nabla \tau }{\partial \tilde{x}_{i}}
\end{equation*}
since 
\begin{equation*}
\int\nolimits_{\gamma _{i}}\frac{\nabla \tau }{\partial \tilde{x}_{i^{\prime
}}}=\frac{\partial \tilde{x}_{i}}{\partial \tilde{x}_{i^{\prime }}}=\delta
_{ii^{\prime }}.
\end{equation*}
Thus under the isomorphism 
\begin{equation}
\left( F^{2}H^{3}\left( \tilde{X}/\tilde{X}^{\prime }\right) \right) ^{\vee }%
\overset{\cong }{\longrightarrow }\Omega _{\tilde{X}^{\prime }}^{1}.
\label{iso2}
\end{equation}
induced by $\left( \ref{Gauss}\right) $ we have 
\begin{equation}
\int\nolimits_{\gamma _{i}}\leftrightarrow d\tilde{x}_{i}\text{.}
\label{corres}
\end{equation}

\subsection{Definition of $\Phi $}

In what follows we will wish to study the function 
\begin{eqnarray}
\Phi &:&\tilde{U}^{\prime }\rightarrow \Bbb{C}  \label{goodfunct} \\
\left( \tilde{u}^{\prime }\right) &\mapsto &\left( \int\nolimits_{r\cdot
\sigma ^{-1}\left( L_{0}\right) \cap X_{u^{\prime }}}^{U_{u^{\prime }}}\tau
\right) .  \notag
\end{eqnarray}
For $\eta $ as above, denote 
\begin{equation}
q=\frac{\tau }{\eta }.  \label{qprime}
\end{equation}
Then $\left( \ref{goodfunct}\right) $ can be rewritten as 
\begin{equation*}
\int\nolimits_{\Gamma }\left( F^{-1}\right) ^{*}\tau ,
\end{equation*}
which, by $\left( \ref{bestyet}\right) $, yields 
\begin{equation}
\Phi =q\int\nolimits_{\Gamma }e^{\left\langle \left. \xi \right| \
\right\rangle }\left( \eta _{0}+\beta \right)  \label{goodform}
\end{equation}
The function $\left( \ref{goodfunct}\right) $ at $\left( 0,0,\eta
_{0}\right) \in \tilde{U}^{\prime }$ can be thought of first as a formal
power series in $\tilde{u}^{\prime }$ and $\overline{\tilde{u}^{\prime }}$,
where as above $\tilde{u}^{\prime }$ is a set of holomorphic coordinates on $%
\tilde{U}^{\prime }$, then as an equivalence class modulo the ideal
generated by the anti-holomorphic coordinates $\overline{\tilde{u}^{\prime }}
$ (see \cite{C}). Notice that 
\begin{equation*}
\tilde{X}^{\prime }
\end{equation*}
is a $\Bbb{C}^{*}$-bundle over $X^{\prime }$. If we denote the (Euler)
vector field which is the derivative of the $\Bbb{C}^{*}$-action as 
\begin{equation*}
\chi ,
\end{equation*}
then by definition, $\Phi $ is a function on the pull-back this bundle to $%
U^{\prime }$ satisfying 
\begin{equation*}
\chi \left( \Phi \right) =\Phi .
\end{equation*}
(That is, $\Phi $ is a section of the dual bundle of $\tilde{U}^{\prime
}/U^{\prime }$.) Let $\Psi $ denote the composite function 
\begin{equation}
F^{2}H^{3}\left( \tilde{X}/\tilde{X}^{\prime }\right) ^{\vee }\rightarrow
F^{3}H^{3}\left( \tilde{X}/\tilde{X}^{\prime }\right) ^{\vee }\overset{%
\left\langle \left. \tau \right| \ \right\rangle }{\longrightarrow }\Bbb{C}
\label{prefunction}
\end{equation}
induced by the inclusion $F^{3}H^{3}\left( \tilde{X}/\tilde{X}^{\prime
}\right) \subseteq F^{2}H^{3}\left( \tilde{X}/\tilde{X}^{\prime }\right) $.
If we let 
\begin{equation}
\varphi =\int\nolimits_{r\cdot \sigma ^{-1}\left( L_{0}\right) \cap
X_{u^{\prime }}}^{U_{u^{\prime }}}\in F^{2}H^{3}\left( X/X^{\prime }\right)
^{\vee }  \label{holmap}
\end{equation}
then 
\begin{equation*}
\Phi =\Psi \circ \varphi .
\end{equation*}
Notice that $\varphi $ is an extension to $U^{\prime }$ of the Abel-Jacobi
map 
\begin{equation*}
\int\nolimits_{r\cdot L\times _{X^{\prime }}Y^{\prime }/Y^{\prime
}}^{Y/Y^{\prime }}
\end{equation*}
on $Y^{\prime }$.

Another way to write the function $\Phi $ is as follows. Use the
trivialization $F$ to rewrite the operator 
\begin{equation}
\varphi =\sum\nolimits_{i=1}^{n^{\prime }}\left( \int\nolimits_{\sigma
^{-1}\Gamma /U^{\prime }}\frac{\nabla \tau }{\partial \tilde{x}_{i}}\right)
\int\nolimits_{\gamma _{i}}=\sum\nolimits_{i=1}^{n^{\prime }}\left(
\int\nolimits_{\Gamma }\left( \left( F^{-1}\right) ^{*}\omega _{i}\right)
\right) \int\nolimits_{\gamma _{i}}.  \label{bettermap}
\end{equation}

Furthermore, referring to $\left( \ref{goodcoords}\right) $ and $\left( \ref
{iso2}\right) $, we define the coordinates 
\begin{equation*}
\left( \tilde{x}_{i}\right) ,\left( w_{i}\right) =\left(
\sum\nolimits_{i=1}^{n^{\prime }}w_{i}d\tilde{x}_{i}\right)
\end{equation*}
for $\left( F^{2}H^{3}\left( \tilde{X}/\tilde{X}^{\prime }\right) \right)
^{\vee }$. In terms of these coordinates, the extension 
\begin{equation*}
\varphi =\int\nolimits_{\sigma ^{-1}\left( \Gamma \right) /U^{\prime
}}=\sum\nolimits_{i=1}^{n^{\prime }}\left( \int\nolimits_{\Gamma }\left(
\left( F^{-1}\right) ^{*}\omega _{i}\right) \right) \int\nolimits_{\gamma
_{i}}
\end{equation*}
of the Abel-Jacobi map on $\tilde{Y}^{\prime }$ is given by 
\begin{equation}
w_{i}=\int\nolimits_{\Gamma }\frac{\partial \left( qe^{\left\langle \left.
\xi \right| \ \right\rangle }\left( \eta _{0}+\beta \right) \right) }{
\partial \tilde{x}_{i}}.  \label{newformula}
\end{equation}

\section{Lagrangian submanifold}

Choose local coordinates $u^{\prime }=\left( \left\{ z_{j}^{\prime }\right\}
,x^{\prime }\right) $ for $U^{\prime }$ where $x^{\prime }$ is the pullbacks
to $U^{\prime }$ of a coordinate system on $X^{\prime }$. If we denote by $%
\tilde{x}_{i}^{\prime }$ the pullbacks to $\tilde{U}^{\prime }$ of the
distinguished coordinate system $\left\{ \tilde{x}_{i}^{\prime }\right\} $
for $\tilde{X}^{\prime }$ defined above, then 
\begin{equation*}
\tilde{u}^{\prime }=\left( \left\{ z_{j}^{\prime }\right\} ,\left\{ \tilde{x}%
_{i}^{\prime }\right\} \right)
\end{equation*}
gives a coordinate system on $\tilde{U}^{\prime }$. Suppose the coordinate
system $\tilde{u}^{\prime }.is$ centered at $\left( 0,0,\eta _{0}\right) $.
We consider the $1$-form 
\begin{equation*}
\Theta \in H^{0}\left( \Omega _{\left( F^{2}H^{3}\left( \tilde{X}/\tilde{X}%
^{\prime }\right) \right) ^{\vee }}^{1}\right)
\end{equation*}
induced by $\left( \ref{iso2}\right) $ and the natural inclusion 
\begin{equation*}
\tilde{\pi }^{*}\Omega _{\tilde{X}^{\prime }}^{1}\subseteq \Omega _{T_{%
\tilde{X}^{\prime }}^{\vee }}^{1}
\end{equation*}
where 
\begin{equation*}
\tilde{\pi }:T_{\tilde{X}^{\prime }}^{\vee }\rightarrow \tilde{X}^{\prime }
\end{equation*}
is the natural projection map on the geometric cotangent bundle. Said
otherwise, under $\left( \ref{iso2}\right) $ we have 
\begin{equation*}
d\tilde{x}_{i}=\int\nolimits_{\gamma _{i}}\nabla \tau \mapsto \left(
\int\nolimits_{\gamma _{i}}:F^{2}H^{3}\left( \tilde{X}/\tilde{X}^{\prime
}\right) \rightarrow \Bbb{C}\right) \in H^{0}\left( \left( F^{2}H^{3}\left( 
\tilde{X}/\tilde{X}^{\prime }\right) \right) ^{\vee }\right) .
\end{equation*}
So the one-form on $T_{\tilde{X}^{\prime }}^{\vee }$ whose value at the
point 
\begin{equation*}
\left( \left( \tilde{x}_{i}\right) ,\left( \sum\nolimits_{i=1}^{n^{\prime
}}w_{i}d\tilde{x}_{i}\right) \right)
\end{equation*}
is 
\begin{equation*}
\sum\nolimits_{i=1}^{n^{\prime }}w_{i}\left( \tilde{\pi }^{*}d\tilde{x}
_{i}\right)
\end{equation*}
corresponds under $\left( \ref{iso2}\right) $ to the one-form $\Theta $
whose value at the point 
\begin{equation*}
\left( \left( \tilde{x}_{i}\right) ,\left( \sum\nolimits_{i=1}^{n^{\prime
}}w_{i}\int\nolimits_{\gamma _{i}}\right) \right)
\end{equation*}
is 
\begin{equation}
\sum\nolimits_{i=1}^{n^{\prime }}w_{i}\left( \tilde{\rho}^{*}d\tilde{x}%
_{i}\right)  \label{form1}
\end{equation}
where 
\begin{equation*}
\tilde{\rho}:\left( F^{2}H^{3}\left( \tilde{X}/\tilde{X}^{\prime }\right)
\right) ^{\vee }\rightarrow \tilde{X}^{\prime }.
\end{equation*}
Then 
\begin{equation*}
d\Theta =\sum\nolimits_{i=1}^{n^{\prime }}dw_{i}\left( \tilde{\rho}^{*}d%
\tilde{x}_{i}\right)
\end{equation*}
descends to give a canonical symplectic structure on the bundle $J\left( 
\tilde{X}/\tilde{X}^{\prime }\right) $ of intermediate Jacobians. In what
follows, we will use the notation $\Theta $ and $d\Psi $ (see $\left( \ref
{prefunction}\right) $ for both the one-forms on $\left( F^{2}H^{3}\left( 
\tilde{X}/\tilde{X}^{\prime }\right) \right) ^{\vee }$ and their pull-backs
to the product $U^{\prime }\times _{X^{\prime }}\left( F^{2}H^{3}\left( 
\tilde{X}/\tilde{X}^{\prime }\right) \right) ^{\vee }$.

\begin{theorem}
\label{goodthm}i) The normal function associated to the family of curves 
\begin{equation*}
\left( Y-r\left( L\times _{X^{\prime }}Y^{\prime }\right) \right) /Y^{\prime
}
\end{equation*}
is the restriction to $\tilde{Y}^{\prime }$ of the holomorphic mapping 
\begin{equation*}
\varphi :\tilde{U}^{\prime }\rightarrow \left( F^{2}H^{3}\left( \tilde{X}/%
\tilde{X}^{\prime }\right) \right) ^{\vee }
\end{equation*}
$\left( \ref{holmap}\right) $. Under the Donagi-Markman isomorphism 
\begin{equation*}
\left( F^{2}H^{3}\left( \tilde{X}/\tilde{X}^{\prime }\right) \right) ^{\vee
}\cong \Omega _{\tilde{X}^{\prime }}^{1},
\end{equation*}
we have the identification 
\begin{equation*}
\varphi \leftrightarrow \varphi ^{*}\Theta .
\end{equation*}

Furthermore, 
\begin{equation*}
\varphi ^{*}\Theta =d_{\tilde{x}^{\prime }}\Phi .
\end{equation*}

ii) The image of $\tilde{U}^{\prime }$ under the ``Abel-Jacobi'' map 
\begin{equation*}
\tilde{U}^{\prime }\overset{\varphi }{\longrightarrow }\left(
F^{2}H^{3}\left( \tilde{X}/\tilde{X}^{\prime }\right) \right) ^{\vee }
\end{equation*}
is ``quasi-Lagrangian,'' that is, 
\begin{equation*}
\varphi ^{*}d\Theta =d_{z^{\prime }}\left( d_{\tilde{x}^{\prime }}\Phi
\right) .
\end{equation*}
\end{theorem}

\begin{proof}
i) The assertion follows directly from $\left( \ref{form1}\right) $ and $%
\left( \ref{newformula}\right) $. Said otherwise, 
\begin{equation*}
\varphi =\sum\nolimits_{i=1}^{n^{\prime }}\left(
\int\nolimits_{rL_{0}}^{Y_{0}}\frac{\partial \left( \left( F^{-1}\right)
^{*}\tau \right) }{\partial \tilde{x}_{i}}\right) \int\nolimits_{\gamma
_{i}},
\end{equation*}
by $\left( \ref{bettermap}\right) $ and, under the isomorphism 
\begin{equation*}
\left( F^{2}H^{3}\left( \tilde{X}/\tilde{X}^{\prime }\right) \right) ^{\vee
}\cong \Omega _{\tilde{X}^{\prime }}^{1},
\end{equation*}
we have, by $\left( \ref{corres}\right) $ that 
\begin{equation*}
\int\nolimits_{\gamma _{i}}\leftrightarrow d\tilde{x}_{i}.
\end{equation*}

ii) 
\begin{eqnarray*}
\varphi ^{*}d\Theta &=&d\varphi ^{*}\Theta \\
&=&d\left( d_{\tilde{X}^{\prime }}\Phi \right) \\
&=&d_{U^{\prime }}\left( d_{\tilde{X}^{\prime }}\Phi \right) .
\end{eqnarray*}
\end{proof}

In our distinguished local coordinates we have by $\left( \ref{newformula}
\right) $ that 
\begin{equation}
\varphi ^{*}\Theta =\sum\nolimits_{i=1}^{n^{\prime }}\left(
\int\nolimits_{\Gamma }\frac{\partial \left( qe^{\left\langle \left. \xi
\right| \ \right\rangle }\left( \eta _{0}+\beta \right) \right) }{\partial 
\tilde{x}_{i}}\right) d\tilde{x}_{i}.  \label{newrformula}
\end{equation}

\section{Hilbert scheme}

\begin{theorem}
\label{zerothm}i) The relative Hilbert scheme $\tilde{Y}^{\prime }$ in $%
\tilde{U}^{\prime }$ is given by the gradient ideal of $\Phi $, that is, the
(formal) germ of $\tilde{Y}^{\prime }$ at $\left( 0,0,\eta _{0}\right) $ is
the zero-scheme of the section 
\begin{equation*}
d_{\tilde{U}^{\prime }/\tilde{X}^{\prime }}\Phi
\end{equation*}
of $\Omega _{\tilde{U}^{\prime }/\tilde{X}^{\prime }}^{1}$.

ii) The image of $\tilde{Y}_{red.}^{\prime }$ under the map 
\begin{equation*}
\varphi :\tilde{U}^{\prime }\rightarrow \left( F^{2}H^{3}\left( \tilde{X}/%
\tilde{X}^{\prime }\right) \right) ^{\vee }
\end{equation*}
is Lagrangian in the sense that the section $\varphi ^{*}\left( d\Theta
\right) $ of $\Omega _{\tilde{U}^{\prime }}^{2}$ restricts to zero on $%
\tilde{Y}_{red.}^{\prime }$. (Compare with \cite{V}, Theorem 1.12.)

iii) For the inclusion 
\begin{equation*}
\iota :\tilde{Y}_{red.sm.}^{\prime }\rightarrow \tilde{U}^{\prime }
\end{equation*}
of the smooth points of $\tilde{Y}_{red.}^{\prime }$, we have 
\begin{equation*}
\left( \varphi \circ \iota \right) ^{*}\left( \Theta \right) =d\left( \iota
^{*}\left( \Phi \right) \right) .
\end{equation*}
\end{theorem}

\begin{proof}
i) Let $\frak{m}$ denote the ideal of $\left\{ 0,0\right\} $ in $U^{\prime }$
. Now 
\begin{equation*}
\left. \xi \right| _{Y_{0}\times X^{\prime }}\in A_{Y_{0}}^{0,1}\left(
T_{Y_{0}}\right) \otimes \Bbb{C}\left[ \left[ u^{\prime }\right] \right] +%
\mathcal{I}_{Y^{\prime }}\cdot A_{Y_{0}}^{0,1}\left( T_{X_{0}}\right)
\end{equation*}
and the obstruction to extending $Y/Y^{\prime }$ to a larger family of
curves is the injective map 
\begin{equation}
\left\{ \xi \right\} :\frac{\mathcal{I}_{Y^{\prime }}}{\frak{m}\cdot 
\mathcal{I}_{Y^{\prime }}}\rightarrow H^{1}\left( N_{Y_{0}\backslash
X_{0}}\right)  \label{ob}
\end{equation}
induced by $\xi $ and the natural map 
\begin{equation*}
\left. T_{X_{0}}\right| _{Y_{0}}\rightarrow N_{Y_{0}\backslash X_{0}}.
\end{equation*}

Under the isomorphism 
\begin{eqnarray*}
T_{\tilde{X}/\tilde{X}^{\prime }} &\rightarrow &\Omega _{\tilde{X}/\tilde{X}
^{\prime }}^{2} \\
\zeta &\mapsto &\left\langle \left. \zeta \right| \tau \right\rangle
\end{eqnarray*}
we can reinterpret $\left( \ref{ob}\right) $ as induced by the map 
\begin{equation}
\frac{\mathcal{I}_{Y^{\prime }}}{\frak{m}\cdot \mathcal{I}_{Y^{\prime }}}%
\rightarrow A^{0,1}\left( \frac{\left. \Omega _{X_{0}}^{2}\right| _{Y_{0}}}{%
\bigwedge\nolimits^{2}N_{Y_{0}\backslash X_{0}}^{\vee }}\right)
=A^{0,1}\left( \left. \Omega _{Y_{0}}^{1}\right| _{Y_{0}}\otimes
N_{Y_{0}\backslash X_{0}}^{\vee }\right)  \label{obs}
\end{equation}
given by 
\begin{equation*}
\left. \left\langle \left. \xi \right| e^{\left\langle \left. \xi \right| \
\right\rangle }\left( \eta _{0}+\beta \right) \right\rangle \right| _{Y_{0}}
\end{equation*}
or, equivalently, by 
\begin{equation*}
\left. \left\langle \left. \xi \right| \eta _{0}\right\rangle \right|
_{Y_{0}}\in \frac{\mathcal{I}_{Y^{\prime }}}{\frak{m}\cdot \mathcal{I}
_{Y^{\prime }}}\cdot A^{0,1}\left( \left. \Omega _{Y_{0}}^{1}\right|
_{Y_{0}}\otimes N_{Y_{0}\backslash X_{0}}^{\vee }\right)
\end{equation*}
since 
\begin{equation*}
\xi ,\beta
\end{equation*}
both vanish along $u^{\prime }=x^{\prime }=0$.

On the other hand, for the function 
\begin{equation*}
\Phi \left( \tilde{u}^{\prime }\right) =\Phi \left( z^{\prime },\tilde{x}
^{\prime }\right) =\int\nolimits_{r\cdot \sigma ^{-1}\left( L_{0}\right)
\cap X_{u^{\prime }}}^{f\left( u^{\prime }\right) }\tau ,
\end{equation*}
when we change $z^{\prime }$ and leave $\tilde{x}^{\prime }$ constant, there
is no change in the complex structure or holomorphic $3$-form $\tau _{\tilde{
x}^{\prime }}$ on $X_{x^{\prime }}$ nor in the position of $\sigma
^{-1}\left( L_{0}\right) \cap X_{x^{\prime }}$. Only the position 
\begin{equation*}
Y_{z^{\prime },x^{\prime }}:=\sigma ^{-1}\left( Y_{0}\right) \cap
X_{x^{\prime }}
\end{equation*}
moves in $X_{x^{\prime }}$. Thus, for a real parameter $t$ we have 
\begin{eqnarray*}
\frac{\partial \Phi }{\partial z_{j}^{\prime }}\left( z^{\prime },\tilde{x}
^{\prime }\right) &=&q\left( \tilde{x}^{\prime }\right) \underset{
t\rightarrow 0}{\lim }\frac{\int\nolimits_{Y_{z^{\prime },x^{\prime
}}}^{Y_{z^{\prime }+t\varepsilon _{j},x^{\prime }}}\eta _{x^{\prime }}}{t} \\
&=&q\left( \tilde{x}^{\prime }\right) \int\nolimits_{Y_{z^{\prime
},x^{\prime }}}\left\langle \left. \gamma _{_{j}\left( z^{\prime },x^{\prime
}\right) }\right| \eta _{x^{\prime }}\right\rangle .
\end{eqnarray*}
Here, as in \S 3 of \cite{C}$,$%
\begin{equation*}
\gamma _{j}\left( z^{\prime },x^{\prime }\right) +\overline{\gamma
_{j}\left( z^{\prime },x^{\prime }\right) }
\end{equation*}
is the real vector field which is the derivative of the family of
diffeomorphisms 
\begin{equation*}
F_{\left( z^{\prime },x^{\prime }\right) }^{-1}\circ F_{\left( z^{\prime
}+t\varepsilon _{j},x^{\prime }\right) }:X_{x^{\prime }}\rightarrow
X_{x^{\prime }}
\end{equation*}
at $t=0$. Recall now that 
\begin{equation*}
\left( F^{-1}\right) ^{*}\left( \eta \right) =e^{\left\langle \left. \xi
\left( z^{\prime },x^{\prime }\right) \right| \ \right\rangle }\left( \eta
_{0}+\beta \left( z^{\prime },x^{\prime }\right) \right)
\end{equation*}
so that 
\begin{equation*}
\frac{\partial \Phi }{\partial z_{j}^{\prime }}\left( z^{\prime },\tilde{x}
^{\prime }\right) =q\left( \tilde{x}^{\prime }\right)
\int\nolimits_{Y_{0}}\left\langle \left. \lambda _{k}\left( z^{\prime
},x^{\prime }\right) \right| e^{\left\langle \left. \xi \left( z^{\prime
},x^{\prime }\right) \right| \ \right\rangle }\left( \eta _{0}+\beta \left(
z^{\prime },x^{\prime }\right) \right) \right\rangle
\end{equation*}
where 
\begin{equation*}
\lambda _{j}\left( z^{\prime },x^{\prime }\right) =\left( F_{\left(
z^{\prime },x^{\prime }\right) }\right) _{*}\left( \gamma _{j\left(
z^{\prime },x^{\prime }\right) }\right)
\end{equation*}
is again of type $\left( 1,0\right) $ by construction. Thus by
considerations of type, we obtain 
\begin{equation}
\frac{\partial \Phi }{\partial z_{j}^{\prime }}\left( z^{\prime },\tilde{x}
^{\prime }\right) =q\left( \tilde{x}^{\prime }\right)
\int\nolimits_{Y_{0}}\left\langle \left. \lambda _{j}\left( z^{\prime
},x^{\prime }\right) \right| \left\langle \left. \xi \left( z^{\prime
},x^{\prime }\right) \right| \left( \eta _{0}+\beta \left( z^{\prime
},x^{\prime }\right) \right) \right\rangle \right\rangle .  \label{important}
\end{equation}

Finally we can consider 
\begin{equation*}
\left\langle \left. \xi \right| e^{\left\langle \left. \xi \right| \
\right\rangle }\left( \eta _{0}+\beta \right) \right\rangle
\end{equation*}
as a section of $\Omega _{U^{\prime }/X^{\prime }}^{1}$ by sending $\frac{
\partial }{\partial z_{k}^{\prime }}$ to 
\begin{equation}
\int\nolimits_{Y_{0}}\left\langle \left. \lambda _{j}\right| \left\langle
\left. \xi \right| \left( \eta _{0}+\beta \right) \right\rangle
\right\rangle .  \label{partials}
\end{equation}
Considering the functions $\left( \ref{partials}\right) $ modulo $\frak{m}
\cdot \mathcal{I}_{\tilde{Y}^{\prime }}$, they can be rewritten as 
\begin{equation*}
\int\nolimits_{Y_{0}}\left\langle \left. \lambda _{j}\left( 0,0\right)
\right| \left\langle \left. \xi \right| \left( \eta _{0}+\beta \right)
\right\rangle \right\rangle .
\end{equation*}

By $\left( \ref{diag}\right) $ the normal vector fields 
\begin{equation*}
\lambda _{j}\left( 0,0\right)
\end{equation*}
give a basis of $H^{0}\left( N_{Y_{0}\backslash X_{0}}\right) $. So we
obtain functions 
\begin{equation}
\left\langle \left. \frac{\partial }{\partial z_{j}^{\prime }}\right| \Phi
\left( y^{\prime },\tilde{x}^{\prime }\right) \right\rangle  \label{gradient}
\end{equation}
which, modulo $\frak{m}\cdot \mathcal{I}_{\tilde{I}^{\prime }}$, are
equivalent to 
\begin{equation*}
z\int\nolimits_{Y_{0}}\left\langle \left. \lambda _{j}\left( 0,0\right)
\right| \left\langle \left. \xi \right| \left( \eta _{0}+\beta \right)
\right\rangle \right\rangle .
\end{equation*}
The injectivity of the map $\left( \ref{ob}\right) $ (rewritten as $\left( 
\ref{obs}\right) $) then implies that the functions 
\begin{equation*}
\int\nolimits_{Y_{0}}\left\langle \left. \lambda _{j}\left( 0,0\right)
\right| \left\langle \left. \xi \right| \left( \eta _{0}+\beta \right)
\right\rangle \right\rangle
\end{equation*}
generate 
\begin{equation*}
\frac{\mathcal{I}_{Y^{\prime }}}{\frak{m}\cdot \mathcal{I}_{Y^{\prime }}}
\end{equation*}
and so, by Nakayama's lemma, generate 
\begin{equation*}
\mathcal{I}_{Y^{\prime }}.
\end{equation*}

ii) By Theorem \ref{goodthm}ii), 
\begin{eqnarray*}
\varphi ^{*}d\Theta &=&d_{z^{\prime }}\left( d_{\tilde{x}^{\prime }}\Phi
\right) \\
&=&-d\left( d_{z^{\prime }}\Phi \right) .
\end{eqnarray*}
But using i) we can conclude that 
\begin{equation*}
\left. \left( image\ of\ d_{z^{\prime }}\Phi \ in\ \Omega _{\tilde{U}
^{\prime }}^{1}\right) \right| _{\tilde{Y}^{\prime }}=0.
\end{equation*}
So the image of $\varphi ^{*}d\Theta $ in $\Omega _{\tilde{Y}
_{red.sm.}^{\prime }}^{2}$ is 
\begin{equation*}
d\left( image\ of\ d_{z^{\prime }}\Phi \ in\ \Omega _{\tilde{Y}
_{red.sm.}^{\prime }}^{1}\right) =0.
\end{equation*}

iii) By Theorem \ref{goodthm} 
\begin{equation*}
\varphi ^{*}\Theta =d_{\tilde{x}^{\prime }}\Phi
\end{equation*}
and by the proof of ii) just above the image of $\iota $ lies in the
zero-locus of $d_{z^{\prime }}\Phi $. So 
\begin{eqnarray*}
\left( \varphi \circ \iota \right) ^{*}\Theta &=&d_{\tilde{x}^{\prime
}}\iota ^{*}\Phi \\
&=&d_{z^{\prime }}\iota ^{*}\Phi +d_{\tilde{x}^{\prime }}\iota ^{*}\Phi \\
&=&\iota ^{*}d\Phi .
\end{eqnarray*}
\end{proof}

\section{The holomorphic Chern-Simons invariant}

The potential function $\Phi $ above is closely related to a variant of the
functional (4) of \cite{DT}. (See also \S 7 of \cite{T}.) To explain this
connection, suggested to the author by R. Thomas, we proceed as follows. Let 
$E_{0}$ be a holomorphic vector bundle on $X_{0}$ (for the moment $X_{0}$
can be any compact complex manifold) and 
\begin{equation*}
\pi :E_{0}^{\vee }\rightarrow X_{0}
\end{equation*}
the dual bundle. For example, $X_{0}$ could be a $K$-trivial threefold and $%
E_{0}$ could be such that, computing its Chern classes as algebraic cycles
modulo rational equivalence, 
\begin{eqnarray}
c_{1}\left( E_{0}\right) &=&0  \label{CSchern} \\
c_{2}\left( E_{0}\right) &\equiv &\left\{ Y_{0}\right\} -r\left\{
L_{0}\right\} .  \notag
\end{eqnarray}
Let 
\begin{equation*}
\rho :U^{\prime }\rightarrow X^{\prime }
\end{equation*}
be a smooth morphism of polydisks as above, but this time let the fiber
dimension of $\rho $ be $h^{1}\left( End\left( E_{0}\right) \right) $, the
embedding dimension at $\left\{ E_{0}\right\} $ of the local analytic scheme
parametrizing holomorphic vector bundles on $X_{0}$. We can compute
deformations of the pair $\left( X_{0},E_{0}\right) $ as follows. Let $%
End^{0}\left( E_{0}\right) $ denote the sheaf of trace-$0$ endomorphisms so
that 
\begin{equation*}
End\left( E_{0}\right) =\mathcal{O}_{X_{0}}\oplus End^{0}\left( E_{0}\right)
.
\end{equation*}
We call the symbol map 
\begin{equation*}
\frak{D}_{1}\left( E_{0}\right) \rightarrow End\left( E_{0}\right) \otimes
T_{X_{0}}
\end{equation*}
followed by the projection 
\begin{equation*}
End\left( E_{0}\right) \otimes T_{X_{0}}\rightarrow End^{0}\left(
E_{0}\right) \otimes T_{X_{0}}
\end{equation*}
the \textit{reduced symbol map}. Call the kernel of the composition 
\begin{equation*}
\frak{D}_{1}^{\dagger }\left( E_{0}\right) .
\end{equation*}
We have the exact sequence of left $\mathcal{O}_{X_{0}}$-modules 
\begin{equation*}
0\rightarrow \frak{gl}\left( E_{0}^{\vee }\right) \rightarrow \frak{D}%
_{1}^{\dagger }\left( E_{0}\right) \rightarrow T_{X_{0}}\rightarrow 0.
\end{equation*}
We consider sections of $E_{0}$ as functions $f$ on $E_{0}^{\vee }$ which
are complex linear on fibers. For any locally defined $C^{\infty }$-vector
field $\beta $ of type $\left( 1,0\right) $ on $X_{0}$ we can lift to a
vector field $\beta ^{\dagger }$ on $E_{0}^{\vee }$ such that 
\begin{equation*}
\pi _{*}\left( \beta ^{\dagger }\right) \in A^{0}\left( \frak{D}%
_{1}^{\dagger }\left( E_{0}\right) \right) .
\end{equation*}

The tangent space to the deformations of $\left( X_{0},E_{0}\right) $ is 
\begin{equation*}
H^{1}\left( \frak{D}_{1}^{\dagger }\left( E_{0}\right) \right)
\end{equation*}
induced by the reduced symbol map. It sits in an exact sequence 
\begin{equation}
H^{1}\left( End\left( E_{0}\right) \right) \rightarrow H^{1}\left( \frak{D}%
_{1}^{\dagger }\left( E_{0}\right) \right) \rightarrow H^{1}\left(
T_{X_{0}}\right) .  \label{CSes}
\end{equation}

Let 
\begin{equation*}
\overline{\partial }_{0}
\end{equation*}
be the $\overline{\partial }$-operator for the complex structure on the
total space of $E_{0}^{\vee }$. The action of $\overline{\partial }_{0}$ on
functions and forms pulled back from $X_{0}$ gives the $\overline{\partial }$%
-operator for the complex structure on $X_{0}$. We fix a $C^{\infty }$
-trivialization 
\begin{equation*}
F:X\rightarrow X^{\prime }\times X_{0}
\end{equation*}
with corresponding Kuranishi data $\xi _{X^{\prime }}$. We embed the
first-order deformations of $\left( X_{0},E_{0}\right) $, which are
parametrized by $H$, into an artinian subscheme of $U^{\prime }$ in a way
compatible with the exact sequence $\left( \ref{CSes}\right) $ and the
differential of $\rho :U^{\prime }\rightarrow X^{\prime }$ at $\left(
0,0\right) $. Next take a maximal extension $E_{Y^{\prime }}/\left(
Y^{\prime }\times _{X^{\prime }}X\right) $ of this first-order deformation
of $\left( X_{0},E_{0}\right) $ over an analytic subscheme $Y^{\prime }$ of $%
U^{\prime }$. Analogous to the deformation of the pair $\left(
Y_{0},X_{0}\right) $ considered earlier, Kuranishi data 
\begin{equation}
\overline{\partial }_{0}-L_{\xi ^{\dagger }}  \label{CSdelbar}
\end{equation}
associated with this maximal formation of $\left( X_{0},E_{0}\right) $
corresponds to a $C^{\infty }$-isomorphism of vector bundles 
\begin{equation}
\begin{array}{ccc}
E:=F^{*}\left( \mathcal{O}_{U^{\prime }}\boxtimes E_{0}\right) & \overset{G}{%
\longrightarrow } & \mathcal{O}_{U^{\prime }}\boxtimes E_{0} \\ 
\downarrow &  & \downarrow \\ 
U^{\prime }\times _{X^{\prime }}X & \overset{\rho ^{*}F}{\longrightarrow } & 
U^{\prime }\times X_{0}
\end{array}
\label{CStriv}
\end{equation}
such that 
\begin{equation*}
E_{Y^{\prime }}\left. :=E\right| _{Y^{\prime }}\times _{X^{\prime }}X
\end{equation*}
has a complex structure induced by the $\overline{\partial }$-operator $%
\left( \ref{CSdelbar}\right) $. Since locally in $X_{0}$, $\xi ^{\dagger }=%
\overline{\partial }_{0}\left( \pi _{*}\beta ^{\dagger }\right) $ for some
vector field $\beta ^{\dagger }$ on $E_{0}^{\vee }$ with $\pi _{*}\beta
^{\dagger }\in A^{0}\left( \frak{D}_{1}^{\dagger }\left( E_{0}\right)
\right) $, $\xi ^{\dagger }$ is a holomorphic mapping 
\begin{equation*}
\xi ^{\dagger }:U^{\prime }\rightarrow A^{0,1}\left( \frak{D}_{1}^{\dagger
}\left( E_{0}\right) \right) \subseteq \pi ^{*}A_{X_{0}}^{0,1}\otimes
T_{E_{0}^{\vee }}
\end{equation*}
lying over (the pull-back to $U^{\prime }$ of) Kuranishi data $\xi $ for the
trivalization $F$ of $X/X^{\prime }$. The scheme $Y^{\prime }$ is then
defined as the maximal solution scheme of the equation 
\begin{equation*}
\left( \overline{\partial }_{0}-L_{\xi ^{\dagger }}\right) ^{2}=0.
\end{equation*}

Next we fix a hermitian metric on $E_{0}$ and, pulling back the hermitian
metric on $E_{0}$ via $\left( \ref{CStriv}\right) $, we induce a hermitian
structure on $E$. Let 
\begin{equation*}
\nabla _{0}:A^{0}\left( E_{0}\right) \rightarrow A^{0}\left( E_{0}\otimes T_{%
\Bbb{C}}^{\vee }\left( X_{0}\right) \right)
\end{equation*}
be the metric $\left( 1,0\right) $-connection for the K\"{a}hler manifold $%
X_{0}$ and the given metric. So, in particular, 
\begin{equation*}
\nabla _{0}^{0,1}=\overline{\partial }_{0}
\end{equation*}
and, with respect to a local holomorphic frame $e_{0}$ of $E,$ the hermitian
structure is given by a hermitian-matrix-valued function $H$ and 
\begin{equation*}
\nabla _{0}^{1,0}\left( e_{0}\right) =\left( \partial H\cdot H^{-1}\right)
\cdot e_{0}=\nabla _{0}\left( e_{0}\right) .
\end{equation*}
Abusing notation, let $\nabla _{0}$ also denote the pull-back of the
connection $\nabla _{0}$ on $E_{0}$ via the product structure on $\left( 
\mathcal{O}_{U^{\prime }}\boxtimes E_{0}\right) /\left( U^{\prime }\times
X_{0}\right) $. By \cite{C} we have 
\begin{eqnarray*}
\nabla _{0}^{1,0} &:&A_{X_{0}}^{0}\left( E_{0}\right) \otimes \Bbb{C}\left[
\left[ U^{\prime }\right] \right] \rightarrow e^{\left\langle \left. \xi
^{\dagger }\right| \ \right\rangle }\left( A_{X_{0}}^{1,0}\left(
E_{0}\right) \otimes \Bbb{C}\left[ \left[ U^{\prime }\right] \right] \right)
\\
\nabla _{0}^{0,1} &=&\overline{\partial }_{0}:A_{X_{0}}^{0}\left(
E_{0}\right) \otimes \Bbb{C}\left[ \left[ U^{\prime }\right] \right]
\rightarrow A_{X_{0}}^{0,1}\left( E_{0}\right) \otimes \Bbb{C}\left[ \left[
U^{\prime }\right] \right] .
\end{eqnarray*}
For a local holomorphic section $s\left( u^{\prime }\right) $ over $%
Y^{\prime }$ we have 
\begin{equation*}
\left( \nabla _{0}-L_{\xi ^{\dagger }}\right) \left( s\left( u^{\prime
}\right) \right) =\nabla _{0}^{1,0}\left( s\left( u^{\prime }\right) \right)
\end{equation*}
so that 
\begin{equation*}
\nabla :=\nabla _{0}-L_{\xi ^{\dagger }}=:d+L_{\alpha }
\end{equation*}
is a $\left( 1,0\right) $-connection. By this we mean 
\begin{equation*}
\nabla ^{0,1}=\overline{\partial }_{0}-L_{\xi ^{\dagger }}.
\end{equation*}

As in Lemma 8.2ii) of \cite{C}, the integrability conditions for the almost
complex structure given by $\left( \overline{\partial }_{0}-L_{\xi ^{\dagger
}}\right) $ on $E_{0}^{\vee }$ are given by the vanishing of the tensor 
\begin{equation*}
\left( \overline{\partial }_{0}-L_{\xi ^{\dagger }}\right) ^{2}.
\end{equation*}
The critical remark is therefore that the $\left( 0,2\right) $ component of
the curvature form $R$ defined by 
\begin{equation*}
L_{R}=\nabla ^{2}=L_{d\alpha +\frac{1}{2}\left[ \alpha ,\alpha \right] }
\end{equation*}
is simply the obstruction tensor 
\begin{equation*}
\left( \overline{\partial }_{0}-L_{\xi ^{\dagger }}\right) ^{2}=-L_{%
\overline{\partial }_{0}\xi ^{\dagger }-\frac{1}{2}\left[ \xi ^{\dagger
},\xi ^{\dagger }\right] }
\end{equation*}
whose vanishing defines $Y^{\prime }$.

Also 
\begin{equation*}
d_{U^{\prime }/X^{\prime }}\left( \overline{\partial }_{0}-L_{\xi ^{\dagger
}}\right) =d\alpha ^{0,1}\in A_{U^{\prime }/X^{\prime }}^{1,0}\otimes
A_{X_{0}}^{0,1}\left( End\left( E_{0}\right) \right)
\end{equation*}
induces the isomorphism 
\begin{equation}
\left. T_{U^{\prime }/X^{\prime }}\right| _{\left( \left\{ X_{0}\right\}
,\left\{ E_{0}\right\} \right) }\rightarrow H^{1}\left( End\left(
E_{0}\right) \right)  \label{CSmap}
\end{equation}
determined by the first-order deformations of $E_{0}$ with $X_{0}$ fixed. If
as above $\left\{ z_{j}^{\prime },x_{i}^{\prime }\circ \rho \right\} $ are a
system of holomorphic coordinates on $U^{\prime }$ we have 
\begin{equation}
\sum\nolimits_{j}\left. \frac{\partial \xi ^{\dagger }}{\partial
z_{j}^{\prime }}\right| _{u^{\prime }=0}dz_{j}^{\prime }  \label{CSb}
\end{equation}
where $\left\{ \left. \frac{\partial \xi ^{\dagger }}{\partial z_{j}^{\prime
}}\right| _{u^{\prime }=0}\right\} $ is a basis of $H^{1}\left( End\left(
E_{0}\right) \right) $.

As before, the maximality of $Y^{\prime }$ implies that the map 
\begin{equation}
\left\{ R^{0,2}\right\} :\left( \frac{\mathcal{I}_{Y^{\prime }}}{\frak{m}
\cdot \mathcal{I}_{Y^{\prime }}}\right) ^{\vee }\rightarrow H^{2}\left(
End\left( E_{0}\right) \right)  \label{CSgee}
\end{equation}
is injective. For a minimal generating set $\left\{ g_{k}\right\} $ of $%
\mathcal{I}_{Y^{\prime }},$ Nakayama's lemma says that the $g_{k}$ are
linearly independent in $\frac{\mathcal{I}_{Y^{\prime }}}{\frak{m}\cdot 
\mathcal{I}_{Y^{\prime }}}$ and so we can write 
\begin{equation}
\left\{ \left( G^{-1}\right) ^{*}R^{0,2}\right\} =\sum\nolimits_{k}\mu
_{k}g_{k}  \label{CSpb}
\end{equation}
for a partial basis $\left\{ \mu _{k}\right\} $ of $H^{2}\left( End\left(
E_{0}\right) \right) $.

Now we return to the situation in which $X_{0}$ is a $K$-trivial threefold.
Following \cite{DT} define the holomorphic Chern-Simons functional as
follows. For each fixed $u^{\prime }=\left( z^{\prime },x^{\prime }\right)
\in U^{\prime }$ embed $\left[ 0,1\right] $ in $U^{\prime }$ via 
\begin{eqnarray}
\left[ 0,1\right] &\rightarrow &U^{\prime }  \label{tmap} \\
t &\mapsto &\left( tz^{\prime },x^{\prime }\right)  \notag
\end{eqnarray}
and form the connection 
\begin{eqnarray*}
\tilde{\nabla}_{u^{\prime }} &:&=\nabla _{\left( tz^{\prime },x^{\prime
}\right) }+d_{t} \\
&=&\nabla _{0}-L_{\xi ^{\dagger }\left( tz^{\prime },x^{\prime }\right)
}+d_{t}
\end{eqnarray*}
on 
\begin{equation*}
\left( \left[ 0,1\right] \cdot z^{\prime },x^{\prime }\right) \times
_{X^{\prime }}X.
\end{equation*}

\begin{definition}
The \textit{holomorphic Chern-Simons functional} $CS_{\nabla _{0}}\left(
u^{\prime }\right) $ on $U^{\prime }$ is given by the formula 
\begin{equation*}
CS_{\nabla _{0}}\left( u^{\prime }\right) =\int\nolimits_{X/X^{\prime }}\tau
\wedge \left( \int\nolimits_{0}^{1}\left( \tilde{R}_{u^{\prime }}\wedge 
\tilde{R}_{u^{\prime }}\right) dt\right) 
\end{equation*}
where 
\begin{equation*}
\tilde{R}_{u^{\prime }}=\tilde{\nabla}_{u^{\prime }}^{2}.
\end{equation*}
\end{definition}

We can decompose $\tilde{R}_{u^{\prime }}$ by type to obtain that 
\begin{equation*}
\left( \left( \nabla _{0}^{1}-L_{\xi ^{\dagger }}+d_{t}\right) ^{2}\right)
^{\left( 0,2\right) ,0+\left( 0,1\right) ,dt}
\end{equation*}
is given by the expression 
\begin{equation*}
R_{0}-\left( \overline{\partial }_{0}\xi ^{\dagger }-\frac{1}{2}\left[ \xi
^{\dagger },\xi ^{\dagger }\right] \right) -\sum\nolimits_{j}z_{j}^{\prime }%
\frac{\partial \xi ^{\dagger }}{\partial z_{j}^{\prime }}dt.
\end{equation*}
So by type 
\begin{eqnarray}
\left( \int\nolimits_{0}^{1}\left( \tilde{R}\wedge \tilde{R}\right)
dt\right) ^{0,3} &=&2\int\nolimits_{0}^{1}\left(
\sum\nolimits_{j}z_{j}^{\prime }\frac{\partial \xi ^{\dagger }}{\partial
z_{j}^{\prime }}\right) \wedge \left( \overline{\partial }_{0}\xi ^{\dagger
}-\frac{1}{2}\left[ \xi ^{\dagger },\xi ^{\dagger }\right] \right) dt
\label{CScomp} \\
&=&2\sum\nolimits_{j}z_{j}^{\prime }\int\nolimits_{0}^{1}\left( \frac{%
\partial \xi ^{\dagger }}{\partial z_{j}^{\prime }}\wedge R^{0,2}\right) dt 
\notag
\end{eqnarray}
where the $\xi ^{\dagger }$ and $R^{0,2}$ under the integral sign are
functions of $t$ via the map $\left( \ref{tmap}\right) $.

\begin{definition}
\label{CSF}We define 
\begin{equation*}
\Phi _{DT}\left( \tilde{u}^{\prime }\right) :=CS_{D_{0}}\left( \nabla ,\tau
\right) 
\end{equation*}
where $D_{0}$ is a ``base'' $\left( 1,0\right) $-connection on $E_{0}$ and
the right-hand expression is the holomorphic Chern-Simons functional
associated to the variation $G^{*}\left( \nabla \right) $ of the connection $%
G^{*}\left( \nabla _{0}\right) $ on the bundle 
\begin{equation*}
G^{*}\left( E_{0}\times U^{\prime }\right) 
\end{equation*}
over 
\begin{equation*}
X\times _{X^{\prime }}U^{\prime }/U^{\prime }.
\end{equation*}
By $\left( \ref{CScomp}\right) $ it is computed as $CS_{D_{0}}\left( \nabla
_{0},\tau \right) +CS_{\nabla _{0}}\left( \nabla ,\tau \right) $ where 
\begin{eqnarray*}
CS_{\nabla _{0}}\left( \nabla ,\tau \right)  &=&\int_{X_{0}\times \left(
\left[ 0,1\right] \cdot u^{\prime }\right) }tr\left( \tilde{R}\wedge \tilde{R%
}\right) \wedge \left( G^{-1}\right) ^{*}\left( \tau \right)  \\
&=&2\int_{X_{0}}\sum\nolimits_{j}z_{j}^{\prime
}\int\nolimits_{0}^{1}tr\left( \frac{\partial \xi ^{\dagger }}{\partial %
z_{j}^{\prime }}\wedge R^{0,2}\left( tu^{\prime }\right) \right) dt\wedge
qe^{\left\langle \left. \xi \right| \ \right\rangle }\left( \eta _{0}+\beta
\right)  \\
&=&2\int_{X_{0}}\sum\nolimits_{j}z_{j}^{\prime
}\int\nolimits_{0}^{1}tr\left( \frac{\partial \xi ^{\dagger }}{\partial %
z_{j}^{\prime }}\wedge R^{0,2}\left( tu^{\prime }\right) \right) dt\wedge
\left( \eta _{0}+\beta \right) .
\end{eqnarray*}
and, as above, $\tau $ is the tautological $\left( 3,0\right) $-form on $%
\tilde{X}/\tilde{X}^{\prime }$ and $\beta $ and $q$ are is as in $\left( \ref
{newrformula}\right) $.
\end{definition}

As before, 
\begin{equation*}
R^{0,2}\left( u^{\prime }\right) \in \frak{m}\cdot \mathcal{I}_{Y^{\prime
}}\cdot A^{0,2}\left( X_{0}\right)
\end{equation*}
where we consider 
\begin{equation*}
\left\{ R^{0,2}\right\}
\end{equation*}
as in $\left( \ref{CSgee}\right) $. And so, modulo $\frak{m}^{2}\cdot 
\mathcal{I}_{Y^{\prime }}$, we have 
\begin{equation*}
\Phi _{DT}\equiv CS_{D_{0}}\left( \nabla _{0},\tau \right)
+q\sum\nolimits_{j}z_{j}^{\prime }\int_{X_{0}}tr\left( \frac{\partial \xi
^{\dagger }}{\partial z_{j}^{\prime }}\left( u^{\prime }\right) \wedge
R^{0,2}\left( u^{\prime }\right) \right) \wedge \eta _{0}.
\end{equation*}
Thus since Serre duality of $H^{1}\left( End\left( E_{0}\right) \right) $
and $H^{2}\left( End\left( E_{0}\right) \right) $ is expressed by the
pairing 
\begin{equation*}
\int_{X_{0}}tr\left\langle \left. {}\right| \ \right\rangle \wedge \eta _{0}
\end{equation*}
we have by $\left( \ref{CSb}\right) $ and $\left( \ref{CSpb}\right) $ we
have modulo $\frak{m}^{2}\cdot \mathcal{I}_{Y^{\prime }}$ that 
\begin{equation*}
\Phi _{DT}\equiv CS_{D_{0}}\left( \nabla _{0},\tau \right)
+q\sum\nolimits_{k}\left( \sum\nolimits_{j}\left\langle \left. \frac{%
\partial \xi ^{\dagger }}{\partial z_{j}^{\prime }}\left( 0\right) \right|
\mu _{k}\right\rangle z_{j}^{\prime }\right) g_{k}.
\end{equation*}
So, modulo $\frak{m}\cdot \mathcal{I}_{Y^{\prime }}$ we have 
\begin{equation*}
\frac{\partial \Phi _{DT}}{\partial z_{j}^{\prime }}=q\sum\nolimits_{j,k}%
\left\langle \left. \frac{\partial \xi ^{\dagger }}{\partial z_{j}^{\prime }}
\left( 0\right) \right| \mu _{k}\right\rangle g_{k}
\end{equation*}
where, by the injectivity of $\left( \ref{CSgee}\right) $, the matrix 
\begin{equation*}
\left\langle \left. \frac{\partial \xi ^{\dagger }}{\partial z_{j}^{\prime }}
\left( 0\right) \right| \mu _{k}\right\rangle
\end{equation*}
is of maximal rank.

Nakayama's lemma then gives an alternate proof of the result, established
independently by Witten and Donaldson-Thomas \cite{DT}, that $\tilde{Y}
^{\prime }$ is the zero-scheme of the section 
\begin{equation*}
d_{\tilde{U}^{\prime }/\tilde{X}^{\prime }}\Phi _{DT}
\end{equation*}
of 
\begin{equation*}
\Omega _{\tilde{U}^{\prime }/\tilde{X}^{\prime }}^{1}.
\end{equation*}
Tjurin's construction in \S 7 of \cite{T} of an Abel-Jacobi map 
\begin{equation*}
Y^{\prime }\cap \rho ^{-1}\left( 0\right) \rightarrow \left(
F^{2}H^{3}\left( X_{0}\right) \right) ^{\vee }
\end{equation*}
is given via $\left( \ref{Gauss}\right) $ by 
\begin{equation*}
\left. d_{\tilde{U}^{\prime }}\Phi _{DT}\right| _{\tilde{Y}^{\prime }\cap 
\tilde{\rho}^{-1}\left( 0\right) }:\tilde{Y}^{\prime }\cap \tilde{\rho}
^{-1}\left( 0\right) \rightarrow \left( F^{2}H^{3}\left( X_{0}\right)
\right) ^{\vee }.
\end{equation*}
Since we are in a relative setting, we can use the isomorphism $\left( \ref
{Gauss}\right) $ to give a more general definition of an Abel-Jacobi or
normal function map. Namely 
\begin{equation*}
\left. d_{\tilde{U}^{\prime }}\Phi _{DT}\right| _{\tilde{Y}^{\prime }}\in 
\tilde{\rho}^{*}\Omega _{\tilde{X}^{\prime }}^{1}
\end{equation*}
and so, via $\left( \ref{Gauss}\right) $, 
\begin{equation*}
\left. d_{\tilde{U}^{\prime }}\Phi _{DT}\right| _{\tilde{Y}^{\prime }}:%
\tilde{Y}^{\prime }\rightarrow \left( F^{2}H^{3}\left( \tilde{X}/\tilde{X}
^{\prime }\right) \right) ^{\vee }.
\end{equation*}

\section{Relation between $\Phi $ and $\Phi _{DT}$, an analogue of Abel's
theorem (by Richard Thomas)}

There is at least one setting in which $\Phi $ and $\Phi _{DT}$ are exactly
the same. Let $E_{0}$ be a $C^{\infty }$ rank-$2$ complex vector bundle with
a fixed trivialisation of its determinant; thus its determinant has a
natural \emph{holomorphic} structure 
\begin{equation*}
\det E_{0}=\mathcal{O}_{X_{0}}.
\end{equation*}
For most of this section we may assume that $X_{0}$ is any compact complex
3-fold. Choosing a hermitian structure on $E_{0}$ reduces the structure
group of $E_{0}$ to $U\left( 2\right) $; the choice of a distinguished 
\begin{equation*}
1\in H^{0}\left( \mathcal{O}_{X_{0}}\right)
\end{equation*}
and \emph{a hermitian metric compatible with this} (i.e. such that $1$ has
constant unit norm in the induced metric on $\Lambda ^{2}E_{0}$) further
reduces the structure group to $SU\left( 2\right) $. Thus $E_{0}$ has the
structure of a (left) quaternionic line bundle with the multiplicative group 
\begin{equation*}
\Bbb{H}^{*}=\left( 0,\infty \right) \times S^{3}
\end{equation*}
of the quaternions 
\begin{equation*}
\Bbb{H}=\left\{ \left[ 
\begin{array}{cc}
a & b \\ 
-\overline{b} & \overline{a}
\end{array}
\right] :a,b\in \Bbb{C}\right\}
\end{equation*}
acting on the left on $E_{0}$. This $\Bbb{H}$-structure is determined by
left multiplication by $\mathbf{j}$ which is given by the equations 
\begin{eqnarray*}
\left( \mathbf{j}\cdot s\right) &\bot &s \\
\left( \mathbf{j}\cdot s\right) \wedge s &=&\left\| s\right\| ^{2}.
\end{eqnarray*}
(Notice that multiplication by $\mathbf{j}=\left[ 
\begin{array}{cc}
0 & -1 \\ 
1 & 0
\end{array}
\right] $ is not $\Bbb{C}$-linear.) The proof of the equality of $\Phi $ and 
$\Phi _{DT}$ is exactly parallel of the classical proof of Abel's theorem in
which $\Bbb{C}$ is replaced by $\Bbb{H}.$ Every $SU(2)$ connection 
\begin{equation*}
D:A_{X_{0}}^{0}\left( E_{0}\right) \rightarrow A_{X_{0}}^{1}\left(
E_{0}\right)
\end{equation*}
is quaternionic, that is, commutes with the left $\Bbb{H}^{*}$-action. Thus
the connection is completely determined by its value on a single
non-vanishing section, just as is the case for $\Bbb{C}$-line bundles. For
any locally defined non-zero $C^{\infty }$-section $s_{0}$ of $E_{0}$ there
is a unique $SU(2)$-connection $D_{0}$ on $E_{0}$ such that 
\begin{equation*}
D_{0}\left( s_{0}\right) =D_{0}^{1,0}\left( s_{0}\right) \oplus
D_{0}^{0,1}\left( s_{0}\right) =0.
\end{equation*}
Then we also have 
\begin{equation*}
D_{0}\left( \mathbf{j}\cdot s_{0}\right) =D_{0}^{1,0}\left( \mathbf{j}\cdot
s_{0}\right) \oplus D_{0}^{0,1}\left( \mathbf{j}\cdot s_{0}\right) =0
\end{equation*}

If $E_{0}$ has a holomorphic structure given by the del-bar operator $%
\overline{\partial }_{0}$ (on sections of $E_{0}$), and $s_{0}$ is a
holomorphic section of $E_{0}$ with zero scheme $Y_{0}$ which is smooth of
codimension $2$, then there is a unique $\Bbb{H}$-connection $D_{0}$ over $%
X_{0}-Y_{0}$ for which $s_{0}$ is flat. $D_{0}$ is not necessarily a $\left(
1,0\right) $-connection. We compute 
\begin{eqnarray*}
s_{0}\wedge \left( \overline{\partial }_{0}-D_{0}^{0,1}\right) \left( 
\mathbf{j}\cdot s_{0}\right) &=&\overline{\partial }_{0}\left( s_{0}\wedge 
\mathbf{j}\cdot s_{0}\right) \\
&=&-\overline{\partial }_{0}\left( \left\| s_{0}\right\| ^{2}\right)
\end{eqnarray*}
so that, in terms of the basis $\left( s_{0},\mathbf{j}\cdot s_{0}\right) $
we have 
\begin{eqnarray}
\overline{\partial }_{0}-D_{0}^{0,1} &=&\left[ 
\begin{array}{cc}
0 & 0 \\ 
\ast & \frac{\overline{\partial }_{0}\left\| s_{0}\right\| ^{2}}{\left\|
s_{0}\right\| ^{2}}
\end{array}
\right]  \notag \\
&=&\left[ 
\begin{array}{cc}
0 & 0 \\ 
\ast & \overline{\partial }_{0}\log \left\| s_{0}\right\| ^{2}
\end{array}
\right] .  \notag
\end{eqnarray}
So for any $\left( 1,0\right) $-connection $\nabla _{0}$, in terms of the
basis $\left( s_{0},\mathbf{j}\cdot s_{0}\right) $ we have 
\begin{eqnarray}
\nabla _{0}^{0,1}-D_{0}^{0,1} &=&\left[ 
\begin{array}{cc}
0 & 0 \\ 
\ast & \overline{\partial }_{0}\log \left\| s_{0}\right\| ^{2}
\end{array}
\right]  \label{CSdiff1} \\
&=&:A^{0,1}.  \notag
\end{eqnarray}
Then for any $\left( 3,0\right) $-form $\tau $, the Chern-Simons functional $%
CS_{\nabla _{0}}\left( D_{0},\tau \right) $ is given by the expression 
\begin{equation}
-\int\nolimits_{X_{0}-Y_{0}}tr\left( A^{0,1}\wedge \overline{\partial }%
_{0}A^{0,1}+\frac{2}{3}\left( A^{0,1}\right) ^{\wedge 3}\right) \wedge \tau
\label{CSdiff"}
\end{equation}
since 
\begin{equation*}
(\nabla _{0}^{2})^{0,2}=\overline{\partial }_{0}^{2}=0.
\end{equation*}

By $\left( \ref{CSdiff1}\right) $, this reduces to 
\begin{equation*}
-\frac{2}{3}\int\nolimits_{X_{0}-Y_{0}}\left( \overline{\partial }_{0}\log
\left\| s_{0}\right\| ^{2}\right) ^{\wedge 3}\wedge \tau =0,
\end{equation*}
and so by the basic additivity $CS_{A}(C)=CS_{A}(B)+CS_{B}(C)$ of
Chern-Simons functionals, 
\begin{equation}
CS_{\nabla _{0}}\left( D_{u^{\prime }},\tau \right) =CS_{D_{0}}\left(
D_{u^{\prime }},\tau \right)   \label{same}
\end{equation}
for any connection $D_{u^{\prime }}$.

Pick a family of sections $\left\{ s_{u^{\prime }}\right\} _{u^{\prime }\in
U^{\prime }}$ of $E_{0}$ (not necessarily holomophic with respect to any
del-bar operator now), with zero set $Y_{u^{\prime }}$. If there happens to
be a finite del-bar operator $\overline{\partial }_{u^{\prime }}$ with
respect to which $s_{u^{\prime }}$ is holomorphic, then these will also
equal the holomorphic Chern-Simons invariant of $\overline{\partial }%
_{u^{\prime }}$, by the above computation. Now $D_{u^{\prime }}-D_{0}$ can
be computed as follows. 
\begin{equation*}
f=s_{u^{\prime }}\cdot s_{0}^{-1}
\end{equation*}
is a $C^{\infty }$-function on $X_{0}-\left( Y_{0}\cup Y_{u^{\prime
}}\right) $ with values in the multiplicative group $\Bbb{H}^{*}$. Then 
\begin{eqnarray*}
0 &=&D_{u^{\prime }}\left( f\cdot s_{0}\right) =df\cdot s_{0}+f\cdot
D_{u^{\prime }}s_{0} \\
D_{u^{\prime }}s_{0} &=&-f^{-1}\cdot df\cdot s_{0}
\end{eqnarray*}
and since $\Bbb{H}$-connections are determined by their values on a single
section 
\begin{eqnarray*}
D_{u^{\prime }}-D_{0} &=&-f^{-1}\cdot df \\
&=&-f^{*}\left( g^{-1}\cdot dg\right)
\end{eqnarray*}
for the invariant one-form 
\begin{equation*}
g=\left[ 
\begin{array}{cc}
a & b \\ 
-\overline{b} & \overline{a}
\end{array}
\right] ^{-1}\cdot d\left[ 
\begin{array}{cc}
a & b \\ 
-\overline{b} & \overline{a}
\end{array}
\right]
\end{equation*}
on $\Bbb{H}^{*}$. Then 
\begin{equation*}
0=d\left( g^{-1}\cdot g\right) =\left( dg^{-1}\right) \cdot g+g^{-1}\cdot dg
\end{equation*}
so that 
\begin{equation*}
d\left( g^{-1}\cdot dg\right) =-g^{-1}\cdot dg\cdot g^{-1}\cdot dg
\end{equation*}
and, \emph{since $D_{0}$ is integrable} $((D_{0}^{2})^{0,2}=0)$, we have the
formula 
\begin{eqnarray}
CS_{D_{0}}\left( D_{u^{\prime }},\tau \right)
&=&\int\nolimits_{X_{0}}tr\left( f^{-1}\cdot df\wedge D_{0}\left(
f^{-1}\cdot df\right) -\frac{2}{3}\left( f^{-1}\cdot df\right) ^{\wedge
3}\right) \wedge \tau  \label{preabel} \\
&=&-\frac{5}{3}\int\nolimits_{X_{0}}f^{*}tr\left( \left( g^{-1}\cdot
dg\right) ^{\wedge 3}\right) \wedge \tau .  \notag
\end{eqnarray}
But 
\begin{equation*}
tr\left( \left( g^{-1}\cdot dg\right) ^{\wedge 3}\right)
\end{equation*}
is the pull-back to $\Bbb{H}^{*}$ of the invariant $3$-form on S$^{3}$ which
generates $H^{3}\left( S^{3}\right) $ and becomes exact on 
\begin{equation*}
S^{3}-\left\{ 1\right\}
\end{equation*}
so that, on $\Bbb{H}^{*}$ we have, for some non-zero constant $c^{\prime }$
that, as distributions, 
\begin{equation*}
tr\left( \left( g^{-1}\cdot dg\right) ^{\wedge 3}\right) \sim c^{\prime
}\int\nolimits_{\left( 0,\infty \right) }.
\end{equation*}
Pulling back via $f$ we therefore have by $\left( \ref{preabel}\right) $
that 
\begin{equation}
CS_{D_{0}}\left( D_{u^{\prime }},\tau \right)
=c\int\nolimits_{L_{0}}^{Y_{u^{\prime }}}\tau .  \label{abel}
\end{equation}
Here we integrate over the 3-chain $f^{-1}(0,\infty )$, which indeed bounds $%
Y_{u^{\prime }}-L_{0}$. This is of course the complete analogue, for rank-$2$
vector bundles on threefolds with trivial determinant, of the classical
Abel's theorem for line bundles on curves.

Suppose now that $X_{0}$ is a Calabi-Yau threefold, and $Y_{0}\subseteq
X_{0} $ is a disjoint unions of smooth elliptic curves, and $Y/Y^{\prime }$
is as above the maximal deformation of $Y_{0}$. Since $H^{1}(\mathcal{O}%
_{X_{0}})=0=H^{2}(\mathcal{O}_{X_{0}})$, Serre's construction gives a unique
holomorphic rank-$2$ vector bundle $E\times _{X^{\prime }}Y^{\prime }$ over $%
X\times _{X^{\prime }}Y^{\prime }$ via the extension 
\begin{equation*}
0\rightarrow \mathcal{O}_{X\times _{X^{\prime }}Y^{\prime }}\rightarrow
E\times _{X^{\prime }}Y^{\prime }\rightarrow \mathcal{I}_{Y,X\times
_{X^{\prime }}Y^{\prime }}\rightarrow 0.
\end{equation*}
Here we assume that 
\begin{equation*}
E\times _{X^{\prime }}Y^{\prime }
\end{equation*}
is the restriction to $X\times _{X^{\prime }}Y^{\prime }$ of some $C^{\infty
}$-vector bundle $E$ on $X\times _{X^{\prime }}U^{\prime }$. Since $%
H^{1}\left( E_{0}\right) =0$, all sections extend and so this family is a
maximal holomorphic deformation of its restriction $E_{0}$ over $X_{0}$.
Then 
\begin{equation*}
\det \left( E\times _{X^{\prime }}Y^{\prime }\right) =\mathcal{O}_{X\times
_{X^{\prime }}Y^{\prime }}
\end{equation*}
and $E\times _{X^{\prime }}Y^{\prime }$ has a section $s$ whose vanishing
scheme is $Y$.

Construct a $C^{\infty }$-trivialization 
\begin{equation*}
F:X\times _{X^{\prime }}U^{\prime }\rightarrow X_{0}\times U^{\prime }
\end{equation*}
such that the maximal deformation $Y/Y^{\prime }$ of $Y_{0}$ satisfies 
\begin{equation*}
Y\subseteq F^{-1}\left( Y_{0}\times U^{\prime }\right)
\end{equation*}
and then construct a compatible trivialization 
\begin{equation*}
G:E\rightarrow E_{0}\times U^{\prime }.
\end{equation*}
We can assume that we have chosen this trivialization so that $G\circ s$ is
given over $Y^{\prime }$ by 
\begin{equation*}
s_{u^{\prime }}=\left( s_{0},u^{\prime }\right) .
\end{equation*}
Then, for the hermitian structure pulled back from $E_{0},$ we have the
hermitian connection $D_{u^{\prime }}$ defined above (such that $%
D_{u^{\prime }}s_{u^{\prime }}=0=D_{u^{\prime }}(\mathbf{j}\cdot
s_{u^{\prime }})$) inducing the possibly singular del-bar operator $%
D_{u^{\prime }}^{0,1}$. By $\left( \ref{same}\right) $ we have $CS_{\nabla
_{0}}(D_{u^{\prime }})=CS_{D_{0}}(D_{u^{\prime }}).$

We therefore have by $\left( \ref{abel}\right) $ that, for $\tilde{u}%
^{\prime }\in \tilde{U}^{\prime }$, 
\begin{equation*}
\Phi _{DT}\left( \tilde{u}^{\prime }\right) =\Phi \left( \tilde{u}^{\prime
}\right) -\Phi \left( 0,0\right) .
\end{equation*}
Differentiating we have 
\begin{equation*}
\left. d\Phi _{DT}\right| _{\tilde{Y}^{\prime }}=\left. d\Phi \right| _{%
\tilde{Y}^{\prime }},
\end{equation*}
that is, that the two Abel-Jacobi maps coincide.

\section{A relative gradient scheme structure on the Noether-Lefschetz locus
(by Claire Voisin)}

Let $X_{0}$ be a Calabi-Yau threefold. Let $L_{0}$ be a very ample line
bundle on $X_{0}$. We are assume that $H^{1}\left( X_{0};\mathcal{O}%
_{X_{0}}\right) =H^{2}\left( X_{0};\mathcal{O}_{X_{0}}\right) =0$ so that
there is no obstruction to deforming $L_{0}$ with $X_{0}$. We are interested
in the deformation theory of the triple $\left( S_{0},X_{0},\lambda \right) $
where 
\begin{equation*}
S_{0}\overset{j}{\hookrightarrow }X_{0}
\end{equation*}
is the inclusion of a smooth member of $\left| L_{0}\right| $, and 
\begin{equation*}
\lambda \in H^{2}\left( S_{0};\Bbb{Z}\right) _{van}\cap H^{1,1}\left(
S_{0}\right)
\end{equation*}
is an integral cohomology class on $S_{0}$ which is both vanishing
(annihilated by $j_{*}$) and of type $\left( 1,1\right) $. Notice that we
have an orthogonal decomposition 
\begin{equation*}
H^{2}\left( S_{0};\Bbb{Q}\right) =\mathrm{image}\left( j^{*}\right) \oplus
H^{2}\left( S_{0},\Bbb{Q}\right) _{van},
\end{equation*}
and that, by the assumption that $H^{2}\left( X_{0};\mathcal{O}
_{X_{0}}\right) =0,$ the first term, namely $\mathrm{image}\left(
j^{*}\right) $, is made of classes which stay of type $\left( 1,1\right) $
under any deformation of the pair $\left( S_{0},X_{0}\right) $. Hence the
restriction to vanishing classes is not really restrictive.

The locus of the pairs $(S_{0},X_{0})$ such that there exists a non zero $%
\lambda $ as above is called the Noether-Lefschetz locus. It splits locally
as a countable union over all $\lambda $'s of the locus where the locally
constant class $\lambda $ remains of $(1,1)$-type. This last locus is called
the component of the Noether-Lefschetz locus determined by $\lambda $.

We make it now more precise. The deformation space of the pair $%
(S_{0},X_{0}) $ is smooth: indeed the deformation theory of $X_{0}$ is
unobstructed, $L_{0} $ deforms with $X_{0}$ (uniquely since $H^{1}\left(
X_{0};\mathcal{O}_{X_{0}}\right) =0,$). And, since by Kodaira vanishing we
have $H^{i}(X_{0};L_{0})=0$,$\,i>0$, the sections of $L_{0}$ deform with $%
X_{0}$. Hence we get a vector bundle over any local deformation space $%
X^{\prime }$ of $X_{0}$, with fiber $H^{0}(X_{x^{\prime }},L_{x^{\prime }})$
over $x^{\prime }\in X^{\prime }$, and the deformations of the pair $%
(S,X_{0})$ are parametrized by a Zariski open set of the projective bundle 
\begin{equation*}
\Bbb{P}H^{0}(X/X^{\prime },L)\rightarrow X^{\prime }
\end{equation*}
with fiber $\Bbb{P}H^{0}(X_{x^{\prime }},L_{x^{\prime }})$.

Let now 
\begin{equation*}
U^{\prime }\subseteq \Bbb{P}H^{0}(X/X^{\prime },L)
\end{equation*}
be an open ball giving the local deformation space of $(S_{0},X_{0})$. As
above we have a smooth morphism 
\begin{eqnarray*}
U^{\prime } &\rightarrow &X^{\prime } \\
u^{\prime } &\mapsto &x^{\prime }=x^{\prime }\left( u^{\prime }\right) .
\end{eqnarray*}
Over $U^{\prime }$, the local systems 
\begin{equation*}
H_{S,\Bbb{Z},van}^{2},\ H_{3,X,S,\Bbb{Z}}
\end{equation*}
with fibers $H^{2}(S_{u^{\prime }},\Bbb{Z})_{van}$ and $H_{3}(X_{x^{\prime
}},S_{u^{\prime }},\Bbb{Z})$ respectively, are trivial. Choose $\lambda \in
H^{2}(S_{0},\Bbb{Z})_{van})\cap H^{1,1}(S_{0})$. The class $\lambda $
provides a locally constant section $(\lambda _{u^{\prime }})_{u^{\prime
}\in U^{\prime }}$ of $H_{S,\Bbb{Z},van}^{2}$ and inside $U^{\prime }$, the
Noether-Lefschetz component determined by $\lambda $ is the scheme 
\begin{equation*}
U_{\lambda }^{\prime }=\{u^{\prime }\in U^{\prime },\,\lambda _{u^{\prime
}}\in F^{1}H^{2}(S_{u^{\prime }})\},
\end{equation*}
where $F^{\cdot }$ denotes the Hodge filtration. This is also equivalent to
the condition 
\begin{equation*}
\lambda _{u^{\prime }}\perp H^{2,0}(S_{u^{\prime }}),
\end{equation*}
where the orthogonality is with respect to the intersection pairing.

Next identifying via Poincar\'{e} duality $H^{2}(S_{u^{\prime }};\Bbb{Z}%
)_{van}$ with 
\begin{equation*}
\ker \left( \,j_{*}\right) :H_{2}(S_{u^{\prime }};\Bbb{Z})\rightarrow
H_{2}(X_{x^{\prime }};\Bbb{Z}),
\end{equation*}
the exact sequence of relative homology provides a surjective map of local
systems 
\begin{equation*}
\partial :H_{3,X,S,\Bbb{Z}}\rightarrow H_{S,\Bbb{Z},van}^{2}.
\end{equation*}
Hence the section $(\lambda _{u^{\prime }})_{u^{\prime }\in U^{\prime }}$
lifts to a section 
\begin{equation*}
(\gamma _{u^{\prime }})_{u^{\prime }\in U^{\prime }}\in H_{3,X,S,\Bbb{Z}}.
\end{equation*}

Note that the homology class which is Poincar\'{e} dual to $\lambda
_{u^{\prime }}$ is represented by a closed differentiable $2$-chain $%
z_{u^{\prime }}$ in $S_{u^{\prime }}$ which is homologous to $0$ in $%
X_{x^{\prime }}$, and that such a lifting $\gamma _{u^{\prime }}$ is
provided by a differentiable $3$-chain $\Gamma _{u^{\prime }}$ in $%
X_{u^{\prime }}$ (varying continuously with $u^{\prime }$) satisfying the
condition 
\begin{equation*}
\partial \Gamma _{u^{\prime }}=z_{u^{\prime }}.
\end{equation*}
As above let 
\begin{equation*}
\tilde{U}^{\prime }\rightarrow \tilde{X}^{\prime }
\end{equation*}
be the deformation spaces of the triples $(S_{0},X_{0},\eta )$ where $\eta $
is a non-zero holomorphic $(3,0)$-form on $X_{0}$. $\tilde{U}^{\prime }$,
resp. $\tilde{X}^{\prime }$, are $\Bbb{C}^{*}$-bundles on $U^{\prime }$,
resp. $X^{\prime }$. Furthermore, as above there is a tautological section $%
\tau $ of the line bundle 
\begin{equation*}
\mathcal{H}^{3,0}\subset \mathcal{H}_{X}^{3},
\end{equation*}
where $\mathcal{H}^{3,0},\,\mathcal{H}_{X}^{3}$ are the Hodge bundles
associated to the family of deformations of $X_{0}$ parametrized by $%
U^{\prime }$ or $X^{\prime }$. Notice that up to this point, everything on $%
U^{\prime }$ is pulled-back from $X^{\prime }$. This is why we use the same
notation. Consider now on $U^{\prime }$ the Hodge bundle $\mathcal{H}
_{X/S}^{3}$ associated with the local system with fiber $H^{3}(X_{u^{\prime
}},S_{u^{\prime }},\Bbb{Z})$. There is a natural bundle map 
\begin{equation*}
f:\mathcal{H}_{X,S}^{3}\rightarrow \mathcal{H}_{X}^{3}
\end{equation*}
on $U^{\prime }$. Now we observe that since each $\tau _{\tilde{u}^{\prime
}}\in H^{3,0}(X_{u^{\prime }})$ vanishes on $S_{u^{\prime }}$, it lifts
naturally to a class in $H^{3}(X_{u^{\prime }},S_{u^{\prime }};\mathbf{C})$.
Hence we get a natural lifting 
\begin{equation*}
\tau _{rel}\in \mathcal{H}_{X,S}^{3}
\end{equation*}
of $\tau $. This section $\tau _{rel}$ is easily seen to be in $F^{3}%
\mathcal{H}_{X,S}^{3}$ where $F^{\cdot }$ here is the Deligne-Hodge
filtration on relative cohomology. Denoting by $<,>$ the pairing between $%
H^{3}(X_{u^{\prime }},S_{u^{\prime }})$ and $H_{3}(X_{u^{\prime
}},S_{u^{\prime }})$, we have now a function $\Phi _{BN}$ on $\tilde{U}
^{\prime }$ which is defined by 
\begin{equation*}
\Phi _{BN}=<\tau _{rel},\tilde{\gamma}>.
\end{equation*}
Here $\tilde{\gamma}$ is the pull-back to $\tilde{U}^{\prime }$ of the
section $(\gamma _{u^{\prime }})_{u^{\prime }\in U^{\prime }}$. Notice that
more concretely, our function can be written as 
\begin{equation*}
\Phi _{BN}(\tilde{u}^{\prime })=\int_{\Gamma _{u^{\prime }}}\tau ,
\end{equation*}
where $\Gamma _{u^{\prime }}$ is a differentiable $3$-chain defined above.

We can now state the following analogue of the results in the previous
sections:

\begin{proposition}
i) The function $\Phi _{BN}$ is holomorphic on $\tilde{U}^{\prime }$.

ii) The inverse image $\tilde{U}_{\lambda }^{\prime }$ of $U_{\lambda
}^{\prime }\subseteq U^{\prime }$ is equal to the relative gradient scheme
associated to the function $\Phi _{BN}$ on $\tilde{U}^{\prime }\rightarrow 
\tilde{X}^{\prime }$, that is 
\begin{equation*}
\tilde{U}_{\lambda }^{\prime }=V(d_{\tilde{U}^{\prime }/\tilde{X}^{\prime
}}\Phi _{BN})
\end{equation*}
where $``V"$ denotes the \textit{vanishing scheme}.
\end{proposition}

\begin{proof}
i) This follows immediately from the formula 
\begin{equation}
\Phi _{BN}=<\tau _{rel},\tilde{\gamma}>,  \label{1.1}
\end{equation}
where $\tilde{\gamma}$ is a flat section of the local system $H_{3,X,S,\Bbb{Z%
}}$ pulled back to $\tilde{U}^{\prime }$, while $\tau _{rel}$ is a
holomorphic section of the bundle $\mathcal{H}_{X,S}^{3}$.

ii) Differentiating the equation (\ref{1.1}), we get 
\begin{equation*}
d_{\tilde{U}^{\prime }/\tilde{X}^{\prime }}\Phi _{BN}=<\nabla _{\tilde{U}%
^{\prime }/\tilde{X}^{\prime }}\tau _{rel},\tilde{\gamma}>,
\end{equation*}
where $\nabla $ is the Gauss-Manin connection on the Hodge bundle $\mathcal{H%
}_{X,S}^{3}$ on $\tilde{U}^{\prime }$. So it suffices to understand 
\begin{equation*}
\nabla _{\tilde{U}^{\prime }/\tilde{X}^{\prime }}\tau _{rel}\in \Omega _{%
\tilde{U}^{\prime }/\tilde{X}^{\prime }}\otimes \mathcal{H}_{X,S}^{3}.
\end{equation*}
The long exact sequence of relative cohomology provides us with a long exact
sequence of Hodge bundles on $\tilde{U}^{\prime },$%
\begin{equation*}
\mathcal{H}_{X}^{2}\overset{j^{*}}{\rightarrow }\mathcal{H}_{S}^{2}\overset{i%
}{\rightarrow }\mathcal{H}_{X,S}^{3}\overset{f}{\rightarrow }\mathcal{H}%
_{X}^{3},
\end{equation*}
where the map $i$ is Poincar\'{e} dual to the map $\partial $ already
considered. We have need the following lemma:

\begin{lemma}
\label{help} 
\begin{equation*}
\nabla _{\tilde{U}^{\prime }/\tilde{X}^{\prime }}\tau _{rel}\in \Omega _{%
\tilde{U}^{\prime }/\tilde{X}^{\prime }}^{1}\otimes i\left( F^{2}\mathcal{H}%
_{S}^{2}\right)
\end{equation*}
Furthermore, the map 
\begin{equation*}
\nabla _{\tilde{U}^{\prime }/\tilde{X}^{\prime }}\tau _{rel}:T_{\tilde{U}%
^{\prime }/\tilde{X}^{\prime }}\rightarrow i\left( F^{2}\mathcal{H}%
_{S}^{2}\right)
\end{equation*}
is surjective.
\end{lemma}

Admitting this lemma, we conclude the proof of the proposition as follows:
We have 
\begin{equation*}
\nabla _{\tilde{U}^{\prime }/\tilde{X}^{\prime }}\tau _{rel}=i(\alpha ),
\end{equation*}
where 
\begin{equation*}
\alpha \in \Omega _{\tilde{U}^{\prime }/\tilde{X}^{\prime }}^{1}\otimes
i\left( F^{2}\mathcal{H}_{S}^{2}\right)
\end{equation*}
is surjective as an element of $Hom\,(T_{\tilde{U}^{\prime }/\tilde{X}
^{\prime }},F^{2}\mathcal{H}_{S}^{2})$. Since $i$ is dual to $\partial $ it
follows that 
\begin{equation*}
<d_{\tilde{U}^{\prime }/\tilde{X}^{\prime }}\tau _{rel},\tilde{\gamma}
>=<\alpha ,\partial \tilde{\gamma}>=<\alpha ,\tilde{\lambda}>,
\end{equation*}
where the pairing in the second and third terms are the intersection pairing
on $\mathcal{H}_{S}^{2}$, and $\tilde{\lambda}$ is the pull-back to $\tilde{U%
}^{\prime }$ of the section $(\lambda _{u^{\prime }})_{u^{\prime }\in
U^{\prime }}$ of $H_{S,\Bbb{Z},van}^{2}$.

It follows that the vanishing of $d_{\tilde{U}^{\prime }/\tilde{X}^{\prime
}}\Phi _{BN}$ at a point $\tilde{u}^{\prime }$ over $u^{\prime }\in
U^{\prime }$ is equivalent to the vanishing of 
\begin{equation*}
<\alpha (\zeta ),\lambda _{u^{\prime }}>,
\end{equation*}
for any $\zeta \in T_{\tilde{U}^{\prime }/\tilde{X}^{\prime },\tilde{u}
^{\prime }}$. But since $\alpha $ is surjective, this is equivalent to the
vanishings 
\begin{equation*}
<\mu ,\lambda _{u^{\prime }}>=0,\,\forall \mu \in H^{2,0}(S_{u^{\prime }}),
\end{equation*}
which provide exactly the equations defining the locus $\tilde{U}_{\lambda
}^{\prime }$.
\end{proof}

We next prove Lemma \ref{help}:

\begin{proof}
By Griffiths transversality, we have 
\begin{equation*}
\nabla _{\tilde{U}^{\prime }/\tilde{X}^{\prime }}\tau _{rel}\in \Omega _{%
\tilde{U}^{\prime }/\tilde{X}^{\prime }}^{1}\otimes F^{2}\mathcal{H}%
_{X,S}^{3}.
\end{equation*}
On the other hand, since $f(\tau _{rel})=\tau $ is pulled-back from $\tilde{X%
}^{\prime }$, the differential $\nabla _{\tilde{U}^{\prime }/\tilde{X}
^{\prime }}\tau _{rel}$ vanishes. Hence we have 
\begin{equation*}
f(\nabla _{\tilde{U}^{\prime }/\tilde{X}^{\prime }}\tau _{rel})=\nabla _{%
\tilde{U}^{\prime }/\tilde{X}^{\prime }}\tau =0,
\end{equation*}
so that 
\begin{equation*}
\nabla _{\tilde{U}^{\prime }/\tilde{X}^{\prime }}\tau _{rel}\in \Omega _{%
\tilde{U}^{\prime }/\tilde{X}^{\prime }}^{1}\otimes i(\mathcal{H}_{S}^{2}).
\end{equation*}
On the other hand, we have 
\begin{equation*}
F^{2}\mathcal{H}_{X,S}^{3}\cap \mathrm{image}\left( \,i\right) =i(F_{S}^{2}%
\mathcal{H}_{S}^{2}).
\end{equation*}
Hence 
\begin{equation*}
\nabla _{\tilde{U}^{\prime }/\tilde{X}^{\prime }}\tau _{rel}\in \Omega _{%
\tilde{U}^{\prime }/\tilde{X}^{\prime }}^{1}\otimes i(F^{2}\mathcal{H}%
_{S}^{2}),
\end{equation*}
which proves the first statement.

It remains to see that 
\begin{equation*}
\nabla _{\tilde{U}^{\prime }/\tilde{X}^{\prime }}\tau _{rel}:T_{\tilde{U}
^{\prime }/\tilde{X}^{\prime }}\rightarrow i(F^{2}\mathcal{H}_{S}^{2})\cong
F^{2}\mathcal{H}_{S}^{2}
\end{equation*}
is surjective. We claim that at a point $\tilde{u}^{\prime }=(u^{\prime
},\tau _{\tilde{u}^{\prime }})\in \tilde{U}^{\prime }$, this maps identifies
up to sign to the composed isomorphism 
\begin{equation*}
H^{0}(S_{u^{\prime }},N_{S_{u^{\prime }}/X_{x^{\prime }}})\overset{\eta
_{u^{\prime }}}{\rightarrow }H^{0}(S_{u^{\prime }},K_{X_{x^{\prime
}}}\otimes N_{S_{u^{\prime }}/X_{x^{\prime }}})\cong H^{0}(S_{u^{\prime
}},K_{S_{u^{\prime }}}).
\end{equation*}
To see this, let $F:=\tilde{U}_{\tilde{x}^{\prime }}^{\prime }$ be the fiber
of the map $\tilde{U}^{\prime }\rightarrow \tilde{X}^{\prime }$ passing
through $\tilde{u}^{\prime }$. Then $F$ identifies to the fiber of $%
U^{\prime }\rightarrow X^{\prime }$ passing through $u^{\prime }$, that is,
to an open set of $\Bbb{P}(H^{0}(X_{x^{\prime }},L_{x^{\prime }}))$, and the
map $\eta _{u^{\prime }}$ is given along $F$ by contraction with a fixed $%
\eta \in H^{3,0}\left( X_{x^{\prime }}\right) $. Let us consider the
universal family 
\begin{equation*}
\mathcal{S}\subset F\times X_{x^{\prime }}.
\end{equation*}
The form $pr_{1}^{*}\eta \in H^{0}(\Omega _{F\times X_{t}}^{3})$ vanishes
under restriction on each fiber of this family, hence provides at $u^{\prime
}$ an element of $\Omega _{F}^{1}\otimes H^{0}(S_{u^{\prime }},\Omega
_{S_{u^{\prime }}}^{2})\cong Hom\,(T_{F},H^{0}(S_{u^{\prime }},\Omega
_{S_{u^{\prime }}}^{2}))$ given by 
\begin{equation}
\zeta \mapsto \left\langle \left. \zeta ^{\prime }\right| \left( \left. {\
pr_{1}^{*}\eta }\right| _{S_{u^{\prime }}}\right) \right\rangle ,
\label{GMmap}
\end{equation}
where $\zeta \in T_{F,u^{\prime }}$ and $\zeta ^{\prime }$ is any $C^{\infty
}$-lifting of $\zeta $ in $\left. {T}_{\mathcal{S}}\right| _{S_{u^{\prime
}}} $. It is a standard fact that up to sign, this element is equal to $%
\nabla _{F}\left( {\eta }\right) $, where $\nabla _{F}$ here is the
Gauss-Manin connection on the bundle of relative cohomology $\mathcal{H}
_{X,S}^{3}$ restricted to $F$.

So we have found that $\nabla _{F}\left( {\eta }\right) \in \Omega
_{F}^{1}\otimes F^{2}\mathcal{H}_{S}^{2}$ is computed at $u^{\prime }$ as
the map $\left( \ref{GMmap}\right) .$ Notice that here $\zeta ^{\prime }$ is
any lifting of $\zeta $ and, in particular, the right hand side can be
computed using local holomorphic liftings of $\zeta $ in ${T}_{\mathcal{S}}$
. But if $\sigma $ is a local defining equation for $\mathcal{S}$, we can
take for $\zeta ^{\prime }$ the vector field 
\begin{equation*}
\zeta ^{\prime }=(\zeta ,0)-\left( 0,\nu \right) ,
\end{equation*}
where $\nu $ is a vertical vector field in $F\times X$, which is normal to $%
\mathcal{S}$ and satisfies $d\sigma (\nu )=d\sigma (\zeta ,0)$. For this
choice we have that 
\begin{equation*}
\left\langle \left. \zeta ^{\prime }\right| \left( \left. {pr_{1}^{*}\eta }%
\right| _{S_{u^{\prime }}}\right) \right\rangle =-\left\langle \left. \nu
\right| \left( \left. {\eta }\right| _{S_{u^{\prime }}}\right) \right\rangle
.
\end{equation*}

Now we note that the map which to $\zeta $ associates $\nu \in
N_{S_{u^{\prime }}/X_{x^{\prime }}}$ such that 
\begin{equation*}
d\sigma (\zeta ,0)=d\sigma (\nu )
\end{equation*}
is exactly the identification of $T_{F,u^{\prime }}$ with $%
H^{0}(S_{u^{\prime }},N_{S_{u^{\prime }}/X_{x^{\prime }}})$, and that the
map which to $\nu $ associates $-\left\langle \left. \nu \right| \left(
\left. {\eta }\right| _{S_{u^{\prime }}}\right) \right\rangle \in
K_{S_{u^{\prime }}}$ is exactly the adjunction isomorphism. Hence our claim
is proved. This concludes the proof of the lemma.
\end{proof}

To complete the parallel with the results of the previous sections, we now
relate $d\Phi _{BN}$ along the Noether-Lefschetz locus with the Abel-Jacobi
map. Note first of all that along the Noether-Lefschetz component $\tilde{U}
_{\lambda }^{\prime }$, the class $\lambda $ is the cohomology class of a
divisor $D_{\lambda }$ on the surface $S_{u^{\prime }}$, which is
well-defined up to rational equivalence, since under our assumptions, the
surfaces $S_{u^{\prime }}$ are regular. This divisor, or $1$-cycle, becomes
homologous to $0$ in $X_{x^{\prime }}$, hence provides a cycle 
\begin{equation*}
Z_{\lambda }=\left( j_{u^{\prime }}\right) _{*}\left( D_{\lambda }\right)
\in CH_{1}(X_{x^{\prime }})_{hom},
\end{equation*}
which has an Abel-Jacobi invariant 
\begin{equation*}
\varphi _{x^{\prime }}(Z_{\lambda })=\int\nolimits_{\Gamma _{u^{\prime
}}}\in \frac{F^{2}H^{3}(X_{x^{\prime }})^{\vee }}{H_{3}(X_{x^{\prime }};\Bbb{%
\ \ Z})}=J(X_{x^{\prime }}).
\end{equation*}
Next, as in $\left( \ref{Gauss}\right) $ there is a natural identification 
\begin{equation}
\nabla \tau :T_{\tilde{X}^{\prime }}\cong F^{2}H^{3}(\tilde{X}/\tilde{X}%
^{\prime })  \label{Gauss'}
\end{equation}
and dually 
\begin{equation*}
\Omega _{\tilde{X}^{\prime }}^{1}\cong F^{2}H^{3}(\tilde{X}/\tilde{X}%
^{\prime })^{\vee },
\end{equation*}
which we can pull-back to $\tilde{U}^{\prime }$.

We have now

\begin{proposition}
At the point $\tilde{u}^{\prime }\in \tilde{U}_{\lambda }$, with images $%
\tilde{x}^{\prime }\in \tilde{X}^{\prime }$ and $u^{\prime }\in U^{\prime }$
, the differential 
\begin{equation*}
\left. d\Phi _{BN}\right| _{\tilde{u}^{\prime }}\in \Omega _{\tilde{X}%
^{\prime },\tilde{x}^{\prime }}^{1}
\end{equation*}
identifies via the above isomorphism to a lifting in $F^{2}H^{3}(X_{x^{%
\prime }})^{\vee }$ of the Abel-Jacobi invariant $\varphi _{x^{\prime
}}(Z_{\lambda })\in J(X_{x^{\prime }})$.
\end{proposition}

\begin{proof}
We have 
\begin{equation*}
d\Phi _{BN}=<\nabla \tau _{rel},\tilde{\gamma}>.
\end{equation*}
But the right hand side is equal to 
\begin{equation*}
\int_{\Gamma _{u^{\prime }}}\left( \nabla \tau _{rel}\right) .
\end{equation*}
We now use the fact that $f(\nabla \tau _{rel})=\nabla \tau $ and the fact
that the integrals 
\begin{equation*}
\int_{\Gamma _{t}}\mu ,\,\mu \in \mathcal{H}_{X,S,u^{\prime }}^{3}
\end{equation*}
vanish on $i(F^{2}\mathcal{H}_{S,u^{\prime }}^{2})=\ker \,f$ to rewrite this
as 
\begin{equation*}
d\Phi _{BN}(\tilde{u}^{\prime })=\int_{\Gamma _{u^{\prime }}}\nabla \tau (%
\tilde{u}^{\prime }).
\end{equation*}
But this means exactly that, via the isomorphism (\ref{Gauss'}) given by $%
\nabla \tau $, the differential $d\Phi _{BN}(\tilde{u}^{\prime })\in \Omega
_{\tilde{X}^{\prime },\tilde{x}^{\prime }}$ identifies with $\int_{\Gamma
_{u^{\prime }}}\in F^{2}H^{3}(Xx^{\prime })^{\vee }$. Since the boundary of
the chain $\Gamma _{u^{\prime }}$ is equal to any differentiable chain
associated to $D_{\lambda }$ by triangulation, the last term is by
definition the Abel-Jacobi invariant of $Z_{\lambda }\in CH_{1}(X_{x^{\prime
}})_{hom}$.
\end{proof}

\end{document}